\documentclass[12pt,a4paper,amstype]{elsart}
\usepackage{amsmath,amssymb,amsfonts,latexsym}
\usepackage{cite}
\usepackage{verbatim}
\usepackage{multicol}
\usepackage{graphicx}
\usepackage{algorithm}
\usepackage{algorithmic}
\usepackage{color}

\newcommand{\ignore}[1]{}

\long\def\comment#1{}

\def \vx{ \boldsymbol{x} }
\def \vp{ \boldsymbol{p} }
\def \vs{ \boldsymbol{s} }

\newcommand{\lon}{\lambda}
\newcommand{\lat}{\theta}
\newcommand{\lond}{\lambda}
\newcommand{\latd}{\theta}
\newcommand{\lons}{\lambda^{\rm e}}
\newcommand{\lats}{\theta^{\rm e}}
\newcommand{\lf}{\left}
\newcommand{\rt}{\right}
\newcommand{\ep}{\varepsilon}
\newcommand{\tngt}{\boldsymbol{\tau}}
\newcommand{\utngt}{\hat{\boldsymbol{\tau}}}
\newcommand{\unrml}{\hat{\boldsymbol{\eta}}}
\newcommand{\frc}{\mathbf{F}}
\newcommand{\matlab}{\textsc{Matlab\;}}

\begin{document}

\begin{frontmatter}

\title{A Study of Different Modeling Choices For Simulating Platelets Within the Immersed Boundary Method}

\author[addr1]{Varun Shankar},
\ead{shankar@cs.utah.edu}
\author[addr2]{Grady B. Wright},
\ead{wright@math.boisestate.edu}
\author[addr3]{Aaron L. Fogelson}
\ead{fogelson@math.utah.edu} and
\author[addr1]{Robert M. Kirby\thanksref{corresp}}
\ead{kirby@cs.utah.edu}

\address[addr1]{School of Computing, Univ.\ of Utah, Salt Lake City, UT, USA}
\address[addr2]{Department of Mathematics, Boise State Univ., Boise, ID, USA}
\address[addr3]{Departments of Mathematics and Bioengineering, Univ.\ of Utah, Salt Lake City, UT, USA}
\thanks[corresp]{Corresponding Author}

\date{July 14, 2011}

\begin{abstract}
  The Immersed Boundary (IB) method is a widely-used numerical
  methodology for the simulation of fluid-structure interaction
  problems.  The IB method utilizes an Eulerian discretization for the
  fluid equations of motion while maintaining a Lagrangian
  representation of structural objects.  Operators are defined for
  transmitting information (forces and velocities) between these two
  representations.  Most IB simulations represent their structures
  with piecewise-linear approximations and utilize Hookean spring
  models to approximate structural forces.  Our specific motivation is
  the modeling of platelets in hemodynamic flows.  In this paper, we
  study two alternative representations -- radial basis functions
  (RBFs) and Fourier-based (trigonometric polynomials and spherical
  harmonics) representations -- for the modeling of platelets in two
  and three dimensions within the IB framework, and compare our
  results with the traditional piecewise-linear approximation
  methodology.  For different representative shapes, we examine the
  geometric modeling errors (position and normal vectors), force computation
  errors, and computational cost and provide an engineering trade-off
  strategy for when and why one might select to employ these different
  representations.
\end{abstract}

\end{frontmatter}
\maketitle

\section{Introduction}
The Immersed Boundary (IB) Method was introduced by Charles Peskin 
in the early 1970's to solve the coupled equations of motion of a
viscous, incompressible fluid and one or more massless, elastic
surfaces or objects immersed in the fluid~\cite{Peskin:1977}.
Rather than generating a body-fitted grid for both exterior and
interior regions of each surface at each timestep and using these
to determine the fluid motion, Peskin instead employed a uniform
Eulerian Cartesian grid over the entire domain and discretized the immersed
boundaries by a set of points that are \textit{not} constrained to
lie on the grid.  In Peskin's work as well as many of the follow-on
works, this set of points was connected via piecewise linear segments
with Hookean spring models being used for approximating structural forces.
Spreading and interpolation operations are then defined for 
transferring force and velocity information between the Lagrangian-defined
structures and the Eulerian-discretized equations of motion.

The IB method was originally developed to model blood flow in the
heart and through heart valves~\cite{Peskin:1977,Peskin:1980,Peskin:1989}, 
but has since been used in a wide variety of other applications, particularly in
biofluid dynamics problems where complex geometries and immersed
elastic membranes or structures are present and make traditional
computational approaches difficult.  Examples include platelet
aggregation in blood clotting~\cite{Fauci_Fogelson:1993,Fogelson:1984,Fogelson2008}, 
swimming of organisms~\cite{Fauci_Fogelson:1993,Fauci_Peskin:1988}, biofilm
processes~\cite{Dillon:1996}, mechanical properties of
cells~\cite{Agresar:1998}, cochlear dynamics~\cite{Beyer:1992},
and insect flight~\cite{Miller_Peskin:2004b,Miller_Peskin:2004}.

We are motivated by the application of the IB method to platelet
aggregation in blood clotting.  Real platelets circulate with
the blood in an inactive state in which they have a discoidal shape.
In order to participate in clot formation, a platelet must undergo
an activation process, one aspect of which is that the platelet
changes shape and becomes more spherical.  In IB modeling, inactive
platelets are approximately elliptical or ellipsoidal in 2D and 3D,
respectively, while activated platelets are approximately circular
in 2D and spherical in 3D.  Piecewise linear approximations
of platelets are currently used within IB methods applied to platelet
aggregation ({\em e.g.}
\cite{Fauci_Fogelson:1993,Fogelson:1984,Fogelson2008}).  We seek
to explore alternative methods for the modeling of platelets that
might decrease the computational time necessary to maintain and update
platelet geometry and motion with comparable or better error
characteristics to the standard piecewise linear models.

In this paper, we examine two alternative representations for
platelets: interpolation with Fourier-based techniques (trigonometric
polynomials in 2D and spherical harmonics in 3D) and interpolation
with radial basis functions (restricted to the unit circle in 2D and
unit sphere in 3D). Fourier methods have frequently been used for the
modeling of circular and spherical objects ({\em
  e.g.}~\cite{mcpeek-shen-farid09a}).  A recent result of Fornberg and
Piret shows that both trigonometric polynomials and spherical
harmonics are just special cases of radial basis functions (RBFs) when
one chooses the shape parameter in a particular limit
\cite{FornbergPiret:2007}.  Additionally, error estimates for RBF
interpolation on the circle and sphere have been given for a much
wider range of target functions than just
$C^{\infty}$~\cite{HubbertMorton:2004,JetterStocklerWard:1999,NarcSunWard:2007}.

To perform a platelet IB computation, one must (1) have a
representation of the surface of the platelet and (2) be able to
compute forces (internal structural forces) at a specified collection
of material points on the platelet surface.  Once forces are
determined, they are ``projected'' to an Eulerian mesh in which they
are incorporated into the solution of the Navier-Stokes equations for
determining the motion of the fluid.  Based upon the updated fluid
velocity field, the platelet's position and shape are updated.  We
will not detail how the projection and interpolation are accomplished
as this has been amply discussed in other works ({\em e.g.}
\cite{Newren2007}).  Our focus is instead restricted to
models for representing the platelet objects and how these can be used
for constructing and maintaining the object's representation,
computing the normal vectors to the object, and computing the internal
structural forces.

For results, we will compare the piecewise linear, Fourier, and RBF
based methods for two different shapes in 2D and two different
shapes in 3D that typify observed platelet geometries.  We
compare the errors in reconstructing these shapes, computing the
normal vectors, and computing the forces.  We provide a discussion of
the engineering trade-offs we observe with respect to error and
computational costs. Our results indicate that the RBF and Fourier
models are viable alternatives to the piecewise linear models for
platelet-like geometries in terms of errors versus computational cost.  We
furthermore find that the RBF models give better results for objects
of varying smoothness than the Fourier models, and thus appear to be
more promising in applications.

The paper is organized as follows.  In Section \ref{sec:modeling} we
present the three different modeling approaches: piecewise linear,
Fourier, and RBFs.  In Section \ref{sec:IBModeling} we review the
components necessary for handling immersed elastic structures in the
IB method.  In Section \ref{sec:implementation} we provide
implementation details for all three models in terms of computing
normal vectors and forces for the platelets.  Results are partitioned
into two sections by spatial dimension.  In Section
\ref{sec:Results2D} we present our comparison of the three modeling
methodologies for 2D platelet objects, while in Section
\ref{sec:Results3D} we present our comparison for 3D platelet objects.
Section \ref{sec:summary} contains a summary of our findings.

\section{Geometric Modeling Strategies}
\label{sec:modeling}

In this section we present the three different geometric modeling
approaches to be examined.  We first present the (traditional)
piecewise linear approach for modeling two and three dimensional
platelet structures.  We then present our two alternative strategies
based on a parametric representation of the surface: Fourier-based
models (trigonometric series in 2D and spherical harmonic series in
3D) and radial basis function (RBF) models. Implementation details for
all three methodologies are provided in Section
\ref{sec:implementation}.


\subsection{Piecewise Linear Models}\label{sec:modeling_pwl}
In the traditional (IB) method, parametric representations of 
the surface are rarely formed explicitly. Typically, a piecewise
linear representation of the boundary is maintained. In 2D, the piecewise linear interpolant is a set of line segments between pairs of IB points.
However, to perform secondary computations (such as computing normals) with a greater level of accuracy than what the piecewise linear interpolant would offer,
piecewise quadratic interpolants are typically fitted to a set of IB points({\em e.g.}~\cite[\S 3.1.1]{Yang2005}). 

Given a parameter \(\lon\), the piecewise quadratic representation is 
therefore defined as:
\begin{align}
x(\lon) &\approx a_x{\lon}^2 + b_x{\lon} + c_x, \label{eq:pquad1} \\
y(\lon) &\approx a_y{\lon}^2 + b_y{\lon} + c_y. \label{eq:pquad2}
\end{align}

The coefficients are computed by solving two linear systems of equations
for each IB point; the right hand sides to these systems of equations are
simply the x and y coordinates of the three IB points to which the piecewise
quadratics are fitted to. Once the coefficients are obtained, one can now
compute derivatives (and therefore normals and other quantities) at each IB
point.

\begin{figure}[htbp]
\begin{center}
\includegraphics[width=0.7\textwidth]{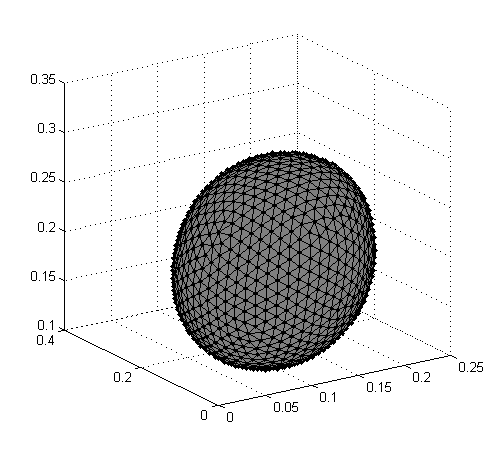}
\end{center}
\caption{Illustration of the triangulation of a set of IB points in 3D.} 
\label{fig:3D_triangles}
\end{figure}

In 3D, the piecewise linear interpolant is a triangulation of the IB points ({\em e.g.}~\cite{Fogelson2008}).  An example of such a
triangulated surface is given in Figure \ref{fig:3D_triangles}.
Secondary computations, such as computing normals and forces are
computed from the triangulation as discussed in Section
\ref{sec:implementation_pwl}.

\subsection{Parametric models}\label{sec:parametric_modeling}
The Fourier and RBF models we propose are both based on explicit
parametric representations of the objects.  Since our target objects
are platelets, which in 2D models are nearly elliptical or
circular and in 3D models are nearly ellipsoidal or spherical, we
choose circular (or polar) and spherical parameterizations in 2D and
3D, respectively.  Before discussing the two modeling approaches, we
introduce some notation and put the modeling problem in the context of
a reconstruction problem using interpolation.

\begin{figure}[htbp]
\begin{center}
\includegraphics[width=0.7\textwidth]{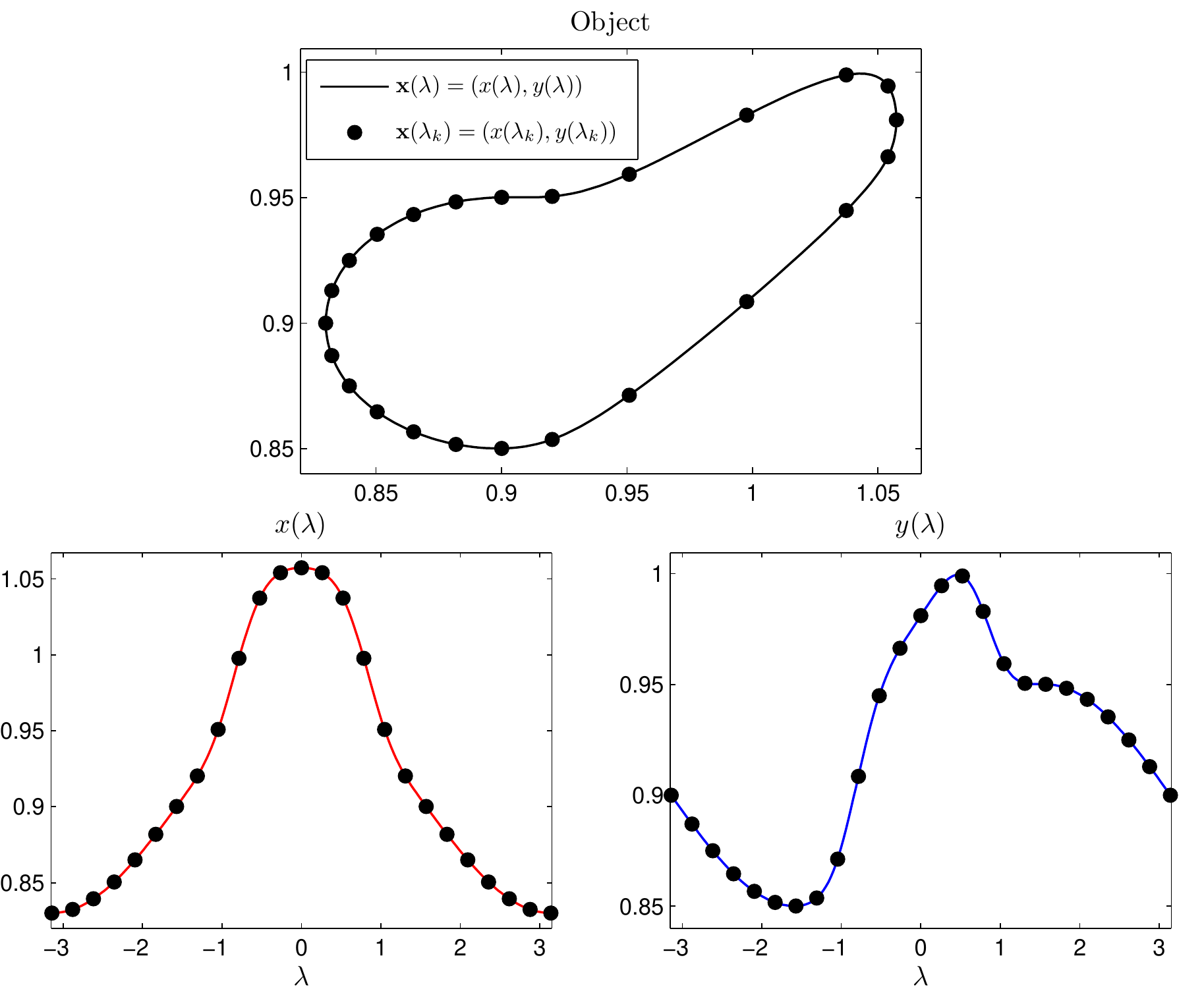}
\end{center}
\caption{Illustration of the parametric representation of a 2D object
  $\vx(\lon)$ and the reconstruction from a finite number of data
  sites.  The top figure shows the 2D object together with discrete data
  sites $\vx(\lond_k) = (x_k,y_k)$.  The bottom left figure shows the $x$
  component of the object in parametric space and its values at the
  node set $\{\lond_k\}_{k=1}^{N}$, while the right figure shows the
  $y$ component and its corresponding values.  The goal is to
  reconstruct $x(\lon)$ and $y(\lon)$ from interpolations of these
  values at the node sets shown and then use these to reconstruct
  $\vx(\lon)$.}
\label{fig:2D_illustration}
\end{figure}

In 2D, we use the following polar parametric notation to represent any of the objects:
\begin{align}
\vx(\lon) = (x(\lon),y(\lon)),
\label{eq:2D_obj}
\end{align}
where $-\pi \leq \lon \leq \pi$ and $\vx(-\pi)=\vx(\pi)$.  In the case
the object is a circle of radius $r$, $\vx(\lon) = (r \cos\lon,r \sin
\lon)$.  In general, given a finite collection of values of the
object, $\{\vx(\lond_k)\}_{k=1}^{N}$ =
$\{(x(\lond_k),y(\lond_k))\}_{k=1}^{N}$, our goal is to reconstruct
$\vx(\lon)$ from smooth interpolations of each of its components.  We
refer to these values as the \emph{data sites} and the set of values
$\{\lond_k\}_{k=1}^N$ as the \emph{nodes}. Figure
\ref{fig:2D_illustration} illustrates this reconstruction problem, of
which the main ingredient is the interpolation of a function defined
on the unit circle.

\begin{figure}[htbp]
\begin{center}
\includegraphics[width=0.85\textwidth]{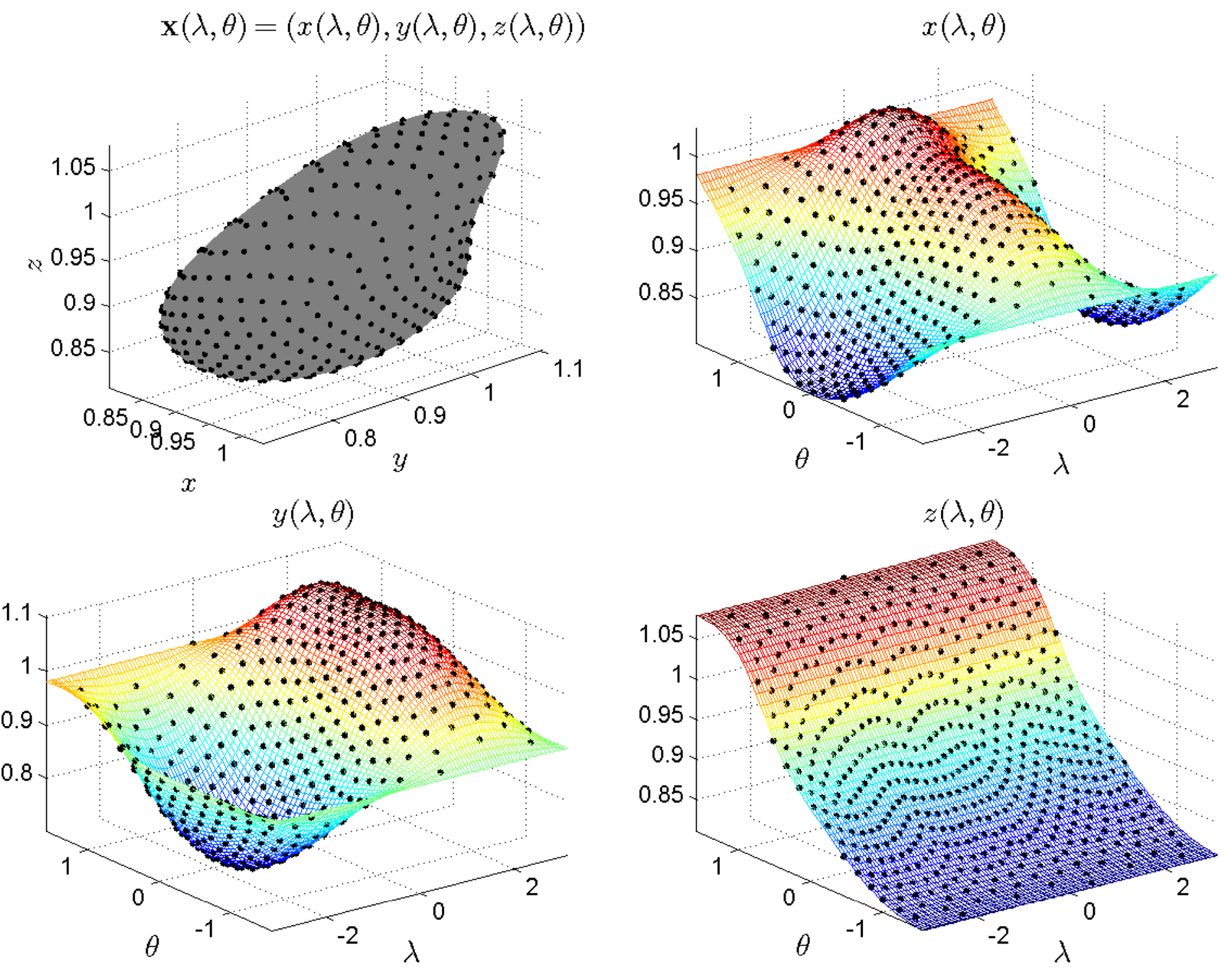}
\end{center}
\caption{Illustration of the parametric representation of a 3D object
  $\vx(\lon,\lat)$ and the reconstruction from a finite number of data
  sites.  Top left figure shows the 3D object together with discrete
  data sites
  $\{(x(\lond_k,\latd_k),y(\lond_k,\latd_k),z(\lond_k,\latd_k)\}_{k=1}^{N}$
  represented as black solid spheres.  Top right figure shows the $x$
  component of the object in spherical parametric space and its values
  at the node set $\{(\lond_k,\latd_k)\}_{k=1}^N$ (marked by black
  solid spheres), while the bottom left and right figures show the
  respective $y$ and $z$ components and their corresponding values.
  The goal is to reconstruct $x(\lon,\lat)$, $y(\lon,\lat)$, and
  $z(\lon,\lat)$ from interpolations of the values at the node sets
  shown and then use these to reconstruct $\vx(\lon,\lat)$.}
\label{fig:3D_illustration}
\end{figure}

In 3D, we represent any of the objects using the following spherical
parametric notation:
\begin{align}
\vx(\lon,\lat) = (x(\lon,\lat),y(\lon,\lat),z(\lon,\lat)),
\label{eq:3D_obj}
\end{align}
where $-\pi \leq \lon \leq \pi$ and $-\pi/2 \leq \lat \leq \pi/2$.
Here the end conditions on $\vx$ in $\lon$ are $\vx(-\pi,\lat) = 
\vx(\pi,\lat)$, while the end conditions in $\lat$ are 
$\vx(\lon,\pi/2)=\vx(\lon+\pi,\pi/2)$ and $\vx(\lon,-\pi/2)=\vx(\lon+\pi,-\pi/2)$
for $-\pi \leq \lon \leq 0$ and $\vx(\lon,\pi/2) = \vx(\lon-\pi,\pi/2)$ and
$\vx(\lon,-\pi/2) = \vx(\lon-\pi,-\pi/2)$ for $0 < \lon \leq \pi$.  These end
conditions on $\theta$ are to enforce continuity of $\vx$ at the poles of the
spherical coordinate system.  In the case the object is a sphere of radius $r$, $\vx(\lon,\lat) =
(r\cos\lon\cos\lat,r\sin\lon\cos\lat,r\sin\lat)$.  Similar to 2D, our
goal is to reconstruct a general object $\vx(\lon,\lat)$ from smooth
interpolations of each of its components which are given at some
finite collection of locations $\{\vx(\lond_k,\latd_k)\}_{k=1}^{N}$ =
$\{(x(\lond_k,\latd_k),y(\lond_k,\latd_k),z(\lond_k,\latd_k))\}_{k=1}^{N}$.
We again refer to these values as the \emph{data sites} and
$\{(\lond_k,\latd_k)\}_{k=1}^N$ as the \emph{nodes}.  Figure
\ref{fig:3D_illustration} illustrates this reconstruction problem, of
which the main ingredient is the interpolation of a function defined
on the unit sphere.

\subsubsection{Fourier Models} \label{sec:modeling_trig} Since the
modeling problems involve interpolation on the unit circle in 2D and
the unit sphere in 3D, a natural choice for constructing these
interpolants are Fourier based methods: trigonometric function for 2D
objects and spherical harmonics for 3D objects.  These methods have
been used extensively for geometric modeling (see for
example~\cite{mcpeek-shen-farid09a} and the references therein).  We
briefly review both of these interpolation techniques in the context
of Figures \ref{fig:2D_illustration} and \ref{fig:3D_illustration}.

Using the notation from Figure \ref{fig:2D_illustration}, we first
discuss the case of reconstructing the $x(\lon)$ component of
$\vx(\lon)$.  In the case that the number of nodes $N$ is even, we
consider a trigonometric interpolant to this data of the form
\begin{align}
p^x(\lon) = c^x_{0} + \sum_{k=1}^{N/2}c^x_{2k-1} \cos k\lon + \sum_{k=1}^{N/2-1}c^x_{2k} \sin k\lon.
\label{eq:trig_interp}
\end{align}
While there is an analogous formula for odd values of $N$, we omit
this discussion and limit our current study to even values of $N$.
The coefficients $c^x_k$ are determined by the interpolation
conditions $p^x(\lond_k) = x(\lond_k)$, $k=1,\ldots,N$.  The solution
to this problem can be written in terms of the following linear
system:

\begin{align}
\begin{bmatrix}
\vspace{0.1in}1 & \cos \lond_1 & \sin \lond_1 & \cdots & \cos \dfrac{N-2}{2}\lond_1 & \sin \dfrac{N-2}{2}\lond_1 & \cos \dfrac{N}{2}\lond_1 \\
1 & \cos \lond_2 & \sin \lond_2 & \cdots & \cos \dfrac{N-2}{2}\lond_2 & \sin \dfrac{N-2}{2}\lond_2 & \cos \dfrac{N}{2}\lond_2 \\
\vdots & \vdots &  \ddots & \vdots & \vdots & \vdots & \vdots \\
1 & \cos \lond_N & \sin \lond_N & \cdots & \cos \dfrac{N-2}{2}\lond_N & \sin \dfrac{N-2}{2}\lond_N & \cos \dfrac{N}{2}\lond_N 
\end{bmatrix}
\begin{bmatrix}
c_0^x\vspace{0.1in} \\
c_1^x\vspace{0.1in} \\
\vdots \\
c_{N-1}^x
\end{bmatrix}
=
\begin{bmatrix}
x_1\vspace{0.1in} \\
x_2\vspace{0.1in} \\
\vdots \\
x_{N}
\end{bmatrix},
\label{eq:trig_linsys}
\end{align}
where $x_k = x(\lond_k)$, $k=1,\ldots,N$. A similar construction to
\eqref{eq:trig_interp} is given for the $y(\lon)$ component of
$\vx(\lon)$, which we denote by $p^{y}(\lon)$.  Our trigonometric
representation of a 2D object like the one in Figure
\ref{fig:2D_illustration} is given by

\begin{align}
{\vp}(\lon) = (p^x(\lon),p^y(\lon)) \label{eq:2D_obj_trig}.
\end{align}

We turn our attention now to interpolation with spherical harmonics
and use the notation from Figure \ref{fig:3D_illustration} to describe
the reconstruction of the $x(\lon,\lat)$ component of
$\vx(\lon,\lat)$.  The dimension of the space of all spherical
harmonics of degree $M$ is given by $(M+1)^2$.  For simplicity, we
thus restrict our attention to the case that the number of nodes is
given by $N = (M+1)^2$.  In this case, we look for a spherical
harmonic interpolant of the form
\begin{align}
p^x(\lon,\lat) = \sum_{\ell=0}^{M}\left[\sum_{m=0}^{\ell}c^x_{\ell,2m} Y_{\ell}^{2m}(\lon,\lat) + \sum_{m=1}^{\ell}c^x_{\ell,2m-1}Y_{\ell}^{2m-1}(\lon,\lat)\right],
\label{eq:sph_interp}
\end{align}
where $Y_{\ell}^{2m}$ and $Y_{\ell}^{2m-1}$ are defined as follows:
\begin{align}
Y_{\ell}^{2m}(\lon,\lat) &:= \sqrt{\frac{2\ell +1}{4\pi}\frac{(\ell-m)!}{(\ell+m)!}}\cos(m\lon)P_{\ell}^{m}(\sin\lat),\quad m=0,\ldots,\ell, \\
Y_{\ell}^{2m-1}(\lon,\lat) &:= \sqrt{\frac{2\ell +1}{4\pi}\frac{(\ell-m)!}{(\ell+m)!}}\sin(m\lon)P_{\ell}^{m}(\sin\lat),\quad m=1,\ldots,\ell.
\end{align}
Here $P_{\ell}^m$ is an associated Legendre function of degree $\ell$
and order $m$.  The coefficients $c^x_k$ are determined by the
interpolation conditions $p^x(\lond_k,\latd_k) = x(\lond_k,\latd_k)$,
$k=1,\ldots,N$.  The linear system corresponding to these conditions
is given by
\begin{align}
\begin{bmatrix}
Y_{0}^{0}(\lond_1,\latd_1) & Y_{1}^{0}(\lond_1,\latd_1) & Y_{1}^{1}(\lond_1,\latd_1) & Y_{1}^{2}(\lond_1,\latd_1) & \cdots \\
Y_{0}^{0}(\lond_2,\latd_2) & Y_{1}^{0}(\lond_2,\latd_2) & Y_{1}^{1}(\lond_2,\latd_2) & Y_{1}^{2}(\lond_2,\latd_2) & \cdots \\
\vdots & \vdots & \vdots & \vdots & \vdots \\
Y_{0}^{0}(\lond_N,\latd_N) & Y_{1}^{0}(\lond_N,\latd_N) & Y_{1}^{1}(\lond_N,\latd_N) & Y_{1}^{2}(\lond_N,\latd_N) & \cdots \\
\end{bmatrix}
\begin{bmatrix}
c^x_1 \\
c^x_2 \\
\vdots \\
c^x_N
\end{bmatrix}
=
\begin{bmatrix}
x_1 \\
x_2 \\
\vdots \\
x_{N}
\end{bmatrix},
\label{eq:sph_linsys}
\end{align}
where $x_k = x(\lond_k,\latd_k)$, $k=1,\ldots,N$. Unlike the
trigonometric case, this linear system can be singular depending on how
the nodes are arranged~\cite[\S 2]{FasshauerSchumaker:1998}.  We avoid
this possibility by choosing the nodes in an ``optimal'' manner as
discussed in Section \ref{sec:implementation_param}.  For a good
review of the properties of spherical harmonic interpolants
see~\cite{WomersleySloan:2001}.  A similar construction to
\eqref{eq:sph_interp} is given for the $y(\lon,\lat)$ and
$z(\lon,\lat)$ components of $\vx(\lon,\lat)$, which we denote by
$p^{y}(\lon,\lat)$ and $p^{z}(\lon,\lat)$.  Our spherical harmonic
representation of a 3D object like the one in Figure
\ref{fig:3D_illustration} is given by

\begin{align}
{\vp}(\lon,\lat) = (p^x(\lon,\lat),p^y(\lon,\lat),p^z(\lon,\lat)).
\label{eq:3D_obj_sph}
\end{align}

\subsubsection{RBF Models}
\label{sec:modeling_rbf}
The RBF method is a popular tool for approximating multidimensional
scattered data. An excellent overview of the theory and application of
this method can be found in the two relatively recent books of
Fasshauer~\cite{Fasshauer:2007} and Wendland~\cite{Wendland:2004}. The
restriction of the RBF method to interpolation on a circle and on a 
sphere began to receive considerable attention from a theoretical
standpoint starting in the mid 1990s (see~\cite[\S
6]{FasshauerSchumaker:1998} for a discussion).  When restricted to
these domains, the RBF method is sometimes referred to as the
\emph{zonal basis function} (ZBF) or \emph{spherical basis function}
(SBF) method in the literature~\cite[Ch. 17]{Wendland:2004}.  We will,
however, use the more popular term RBF to describe the interpolation
technique.  Several studies have been devoted to providing error
estimates for RBF interpolation on circles and spheres; see, for
example,~\cite{JetterStocklerWard:1999,NarcSunWard:2007}.  In the
first of these papers, it is shown these interpolants can provide
spectral accuracy provided the underlying target function is
sufficiently smooth.  The latter of these studies gives error
estimates in the case that the target function belongs to some Sobolev
space.  Recently, the RBF method has been successfully used for
approximating derivatives of scalar and vector-valued quantities on
the surface of a sphere and incorporated into methods for
solving partial differential equations numerically in
spherical
geometries~\cite{Gia:2005,FlyerWright:2007,FlyerWright:2009}.

The construction of the 2D and 3D RBF models of the objects is similar, so we discuss them together.  Using the notation of Figures \ref{fig:2D_illustration} and \ref{fig:3D_illustration}, and focusing on the reconstructions of the $x(\lon)$ and $x(\lon,\lat)$ components of the objects, the corresponding RBF interpolants are given by
\begin{eqnarray}
\text{2D}:&\quad s^x(\lon) &= \sum_{k=1}^{N} c^x_k \phi\lf(\sqrt{2 - 2\cos(\lon-\lond_k)}\rt),
\label{eq:circ_rbf_interp} \\
\text{3D}:&\quad s^{x}(\lon,\lat) &= \sum_{k=1}^{N} c^{x}_k \phi\lf(\sqrt{2(1-\cos \lat \cos \latd_{k}\cos (\lon -\lond_{k})-\sin\lat \sin \latd_{k})}\rt).
\label{eq:sph_rbf_interp}
\end{eqnarray}
Here $\phi$ is some scalar-valued, positive (semi-) definite radial
kernel.  The square root term in \eqref{eq:circ_rbf_interp} is just
the Euclidean distance between the points described in polar
coordinates by $\lon$ and $\lond_k$, while the square root term in
\eqref{eq:sph_rbf_interp} is similarly the Euclidean distance between
the points described in spherical coordinates by $(\lon,\lat)$ and
$(\lond_k,\latd_k)$.  The coefficients $c^x_k$ in either
\eqref{eq:circ_rbf_interp} or \eqref{eq:sph_rbf_interp} are again
determined by the interpolation conditions.  These conditions lead to
the following linear system of equations:
\begin{align}
\underbrace{
\begin{bmatrix}
\phi\lf(r_{1,1}\rt) & \cdots & \phi\lf(r_{1,N}\rt) \\
\phi\lf(r_{2,1}\rt) & \cdots & \phi\lf(r_{2,N}\rt) \\
\vdots & \ddots & \vdots \\
\phi\lf(r_{N,1}\rt) & \cdots & \phi\lf(r_{N,N}\rt)
\end{bmatrix}}_{\displaystyle A}
\begin{bmatrix}
c^x_1 \\
c^x_2 \\
\vdots \\
c^x_N
\end{bmatrix}
= 
\begin{bmatrix}
x_1 \\
x_2 \\
\vdots \\
x_N
\end{bmatrix},
\label{eq:rbf_linsys}
\end{align}
where $x_k = x(\lond_k)$ and $r_{j,k} = \sqrt{2 - 2\cos(\lond_j -
  \lond_k)}$ for 2D objects, and $x_k = x(\lond_k,\latd_k)$ and
$r_{j,k}$ $= \sqrt{2(1-\cos \latd_j \cos \latd_{k}\cos (\lond_j
  -\lond_{k})-\sin\latd_j \sin \latd_{k})}$ for 3D objects.  Note that
$r_{j,k} = r_{k,j}$ so that the linear system \eqref{eq:rbf_linsys} is
symmetric.  More importantly, this linear system is guaranteed to be
non-singular for the appropriate choice of $\phi$.  In this study we
restrict our attention to the multiquadric (MQ) and inverse
multiquadric (IMQ) radial kernels, which are popular in applications
and are given explicitly by
\begin{align}
\text{MQ:}\quad & \phi(r) = \sqrt{1 + (\ep r)^2}, \label{eq:mq} \\
\text{IMQ:}\quad& \phi(r) = \dfrac{1}{\sqrt{1 + (\ep r)^2}}. \label{eq:imq}
\end{align}
Here $\ep$ is called the shape parameter.  For both the MQ and IMQ,
the linear system \eqref{eq:rbf_linsys} is guaranteed to be
non-singular (provided $\ep > 0$).  Furthermore, for the IMQ, the $A$
matrix in this linear system is guaranteed to be positive definite.  A
full discussion of the non-singularity of \eqref{eq:rbf_linsys} for
various radial kernels can be found in either~\cite{Fasshauer:2007}
or~\cite{Wendland:2004}.  We postpone the discussion of choosing $\ep$
to Section \ref{sec:implementation_rbf}.  We do, however, note that in
the limit that $\ep\rightarrow 0$ a RBF interpolant on a circle
converges to a trigonometric interpolant, while a RBF interpolant on a
sphere converges to a spherical harmonic
interpolant~\cite{FornbergPiret:2007} (strictly speaking this was only
shown for the case of the sphere, but the arguments
from~\cite{FornbergPiret:2007} carry directly over to the case of the
circle as well).  Thus, trigonometric and spherical harmonic
interpolation can be viewed as a special case of RBF interpolation.

We denote the RBF representations of a 2D object like the one in
Figure \ref{fig:2D_illustration} by

\begin{align}
{\vs}(\lon) = (s^x(\lon),s^y(\lon)),
\label{eq:2D_obj_rbf}
\end{align}

where $s^y(\lon)$ interpolates $y(\lon)$ and has the form of \eqref{eq:circ_rbf_interp}.  Similarly we denote the RBF representation of a 3D object like the one in Figure \ref{fig:3D_illustration} by

\begin{align}
{\vs}(\lon,\lat) = (s^x(\lon,\lat),s^y(\lon,\lat),s^z(\lon,\lat)),
\label{eq:3D_obj_rbf}
\end{align}
where $s^y(\lon,\lat)$ and $s^z(\lon,\lat)$ interpolate $y(\lon,\lat)$ and $z(\lon,\lat)$, respectively, and have the form of \eqref{eq:sph_rbf_interp}.  

We conclude this section by noting that the RBF method is more
flexible than the Fourier-based methods in regard to altering the
parameterization for the objects.  For example, if one were to find
that a more general ellipse or ellipsoid provided a better
parameterization of the object than a circle or sphere, then the RBF
method can be naturally extended to this new parameterization.  The
only change to \eqref{eq:circ_rbf_interp} or \eqref{eq:sph_rbf_interp}
would be to replace the distance measure in the argument of
$\phi$ with the appropriate (Euclidean) distance measure on the target
object for the parametrization.  More general objects, including ones
with higher genus, are also possible; see~\cite{FuselierWright:2010}
for a theoretical and numerical discussion.

\section{Immersed Boundary Modeling}
\label{sec:IBModeling}
In this section we review the components needed for an immersed
boundary model of platelets. Our focus here is on the computation of
normal vectors and on the modeling of elasticity.
For a discussion of how
forces generated from immersed objects are transferred to the
underlying Eulerian mesh and how the fluid velocity is updated see,
for example,\cite{Newren2007}.

Normal vectors do not play a prominent role in most
traditional IB calculations. Our interest in them is
motivated by their other uses in the modeling of platelet aggregation.
In addition to the fluid-structure interactions modeled using the IB
method, the platelet problem requires solution of
advection-diffusion equations for chemicals in the fluid domain
outside of the moving platelets, along with boundary conditions on
the chemical concentration at the fluid-platelet interface.  Normal
vectors along the platelet boundary are needed for determining when
an Eulerian grid point is inside or outside of the platelet, and for
imposing the boundary conditions.  For further discussion of this,
see \cite{YaoFogelson2011}.

\subsection{Components for 2D}
We denote the 2D platelet using the parametric representation
$\vx(\lon)$ given in \eqref{eq:2D_obj} and define
\begin{align}
\tngt := \frac{\partial}{\partial \lon} \vx(\lon) = \lf(\frac{\partial}{\partial \lon} x(\lon),\frac{\partial}{\partial \lon} y(\lon)\rt) = (\tngt_x,\tngt_y).
\label{eq:tngt_2d}
\end{align}
The unit tangent and normal vectors to $\vx(\lon)$ are then given as
\begin{align}
\utngt: &= \frac{\tngt}{\|\tngt\|} = (\utngt_x,\utngt_y),  \label{eq:utngt_2d} \\
\unrml: &= (-\utngt_y,\utngt_x) \label{eq:unnrml_2d}
\end{align}
For the force model in 2D, we use the fiber model defined
in~\cite{Peskin:2002}.  According to this model, the elastic force
density on $\vx(\lon)$ at the location $\vx(\lon_i)$ is given by
\begin{align}
\frc(\vx(\lon_i)) = \left.\frac{\partial}{\partial \lon} (T\utngt)\right|_{\lon_i}, \label{eq:force_full_2d}
\end{align}
\noindent where $T = K (\|\tngt\|)$ is the fiber tension.  In our
platelet model, we choose $K$ as a linear function, $K =
K_0\|\tngt\|$, where $K_0$ is the Hookean spring constant.  In this
case, \eqref{eq:force_full_2d} reduces to
\begin{align}
\frc(\vx(\lon_i)) = K_0 \left.\frac{\partial}{\partial \lon} \lf(\|\tngt\| \utngt\rt) \right|_{\lon_i} = K_0 \left.\frac{\partial^2}{\partial \lon^2} \vx(\lon) \right|_{\lon_i}. \label{eq:force_2d}
\end{align}

The 2D spring force model traditionally used in piecewise linear
representations is a scaled second-order, central-difference
approximation to the above fiber model (assuming springs of zero rest
length). From the physical standpoint, each IB point in a 2D object is
thought to be connected to each of its neighbors via springs. For
tension forces, there are only two neighbors attached to each IB point
via springs. This spring force is expressed as:

\begin{equation}
\frc(\vx_i) = K_0(\vx_{i+1} - 2\vx_i + \vx_{i-1}). \label{eq:spring_2d}
\end{equation}

\subsection{Components for 3D}
We denote the 3D platelet using the parametric representation $\vx(\lon,\lat)$ given in \eqref{eq:3D_obj} and define
\begin{align}
\tngt^{\lon} := & \frac{\partial}{\partial \lon} \vx(\lon,\lat) = \lf(\frac{\partial}{\partial \lon} x(\lon,\lat),\frac{\partial}{\partial \lon} y(\lon,\lat),\frac{\partial}{\partial \lon} z(\lon,\lat)\rt), \label{eq:tngt1_3D}\\
\tngt^{\lat} := & \frac{\partial}{\partial \lat} \vx(\lon,\lat) = \lf(\frac{\partial}{\partial \lat} x(\lon,\lat),\frac{\partial}{\partial \lat} y(\lon,\lat),\frac{\partial}{\partial \lat} z(\lon,\lat)\rt). \label{eq:tngt2_3D}
\end{align}
The unit tangent vectors to $\vx(\lon,\lat)$ are then given by
\begin{align}
\utngt^{\lon} := \frac{\tngt^{\lon}}{\|\tngt^{\lon}\|}\quad\text{and}\quad\utngt^{\lat} := \frac{\tngt^{\lat}}{\|\tngt^{\lat}\|}, \label{eq:utngt_3d}
\end{align}
while the unit normal vector is given by
\begin{align}
\unrml := \frac{\tngt^{\lon} \times \tngt^{\lat}}{\|\tngt^{\lon} \times \tngt^{\lat}\|}.
\label{eq:unrml_3d}
\end{align}

The force model we use in 3D differs depending on whether a piecewise
linear representation for the object is used or a parametric
representation. Traditionally, piecewise linear representations
(triangulated surfaces) in 3D have been used in conjunction with
spring forces. In this model, a spring is assumed to be placed along
each triangle edge (again, we assume these springs have a rest length
of zero). Then, the total force acting on an IB point at $\vx_i$ due
to its \(k\) nearest neighbors is:
\begin{equation}
\frc(\vx_i) = K_0\sum_{j\neq i}(\vx_i - \vx_j), \label{eq:spring_3d}
\end{equation}
where the sum is over $k$ IB points. The nearest neighbors are
typically defined from the triangulation, i.e. as members of the
adjacency list of $\vx_i$.  This is the same strategy that we follow.

For our parametric representation of platelets, we use surface tension
as the model to compute tension forces.  
The force due to surface tension is then given by
\begin{align}
\frc = \gamma(2H)\unrml,  \label{eq:force_3d}
\end{align}
where $\gamma$ is the coefficient of surface tension.  $H$ is the mean
curvature of the surface, and can be computed as~\cite[\S 16.5]{Gray:1997}
\begin{align}
H = \frac{eG-2fF+gE}{2(EG-F^2)}, \label{eq:mean_curvature}
\end{align}
where $E$, $F$, and $G$ are coefficients of the first fundamental form,
\begin{align}
E = \tngt^{\lon} \cdot \tngt^{\lon},\; F = \tngt^{\lon} \cdot \tngt^{\lat},\; G = \tngt^{\lat} \cdot \tngt^{\lat},
\label{eq:first_fund_form}
\end{align}
and $e$, $f$, and $g$ are coefficients of the second fundamental form,
\begin{align}
e = \lf(\frac{\partial}{\partial \lon}\tngt^{\lon} \rt)\cdot \unrml,\; f = \lf(\frac{\partial}{\partial \lat}\tngt^{\lon} \rt)\cdot \unrml,\;
g = \lf(\frac{\partial}{\partial \lat}\tngt^{\lat} \rt)\cdot \unrml. \\
\label{eq:second_fund_form}
\end{align}
\section{Implementation Details}
\label{sec:implementation}

In this section we present the implementation details for evaluating 
the positions on the Lagrangian objects,  computing normals to the surface 
of the object, and computing the internal forces as presented in the previous 
section.  For the piecewise linear representation, these surface normals and 
forces are computed at the IB points.  For the parametric 
representations using Fourier and RBF models, these values are computed at 
some set of \emph{sample sites}, which do not necessarily correspond to 
the \emph{data sites}. With these operations defined, it is possible to 
employ the traditional spreading and interpolation operators for transferring 
the forces and velocity respectively between the Lagrangian and Eulerian 
discretizations.

\subsection{Piecewise Linear Models}
\label{sec:implementation_pwl}

In 2D, normals are computed at the IB points using the piecewise
quadratic representation presented in Section
\ref{sec:modeling_pwl}. For each IB point, we first solve for the
coefficients in \eqref{eq:pquad1} and \eqref{eq:pquad2} using the IB
point and its two neighbors. Using \eqref{eq:pquad1} and
\eqref{eq:pquad2}, we next compute the tangent vector at each IB point
using \eqref{eq:utngt_2d} and then determine the normal vector using
\eqref{eq:unnrml_2d}.

In 3D, we compute the normal vectors at each IB point by first
computing the normal vector at the circumcenter of each of
the triangles.  We then obtain the normal vector at a vertex (IB point) by
a weighted average of the values of the normal vectors 
at the circumcenters of the triangles connected to the vertex. Specifically,
we weight these facet normals by the angle at which that facet is incident
on the vertex at which we require a normal. This approach takes into account
the geometric configuration of each facet~\cite{ThurmerWuthrich98}.

The implementation of the forces follows directly from the simple
spring force model in both 2D \eqref{eq:spring_2d} and 3D
\eqref{eq:spring_3d}. We note that while the 2D implementation follows
naturally from a constitutive model, the 3D implementation is a purely
algorithmic extension of the 2D case. 

\subsection{Parametric Models}
\label{sec:implementation_param}
For the parametric models, we use the continuous representations of
the objects from either the Fourier or RBF based interpolants to
approximate the normal vectors and forces.  This involves analytically
computing derivatives of these interpolants and then evaluating the
derivatives at some set of $M$ locations in the parametric space that
corresponds to the set of sample sites. In 2D, we denote the set of
sample sites by $\{\vx(\lons_j)\}_{j=1}^M$ and refer to the set of
parametric values $\{\lons_j\}_{j=1}^M$ as the \emph{evaluation
  points}.  Similarly for 3D, we denote the sample sites by
$\{\vx(\lons_j,\lats_j)\}_{j=1}^M$ and refer to
$\{(\lons_j,\lats_j)\}_{j=1}^M$ as the \emph{evaluation points}. The
method we use is similar to the pseudospectral or spectral collocation
method ({\em e.g.}~\cite{Fornberg:1996:PGPM,Trefethen:2000:SMM}),
except that the derivatives are not evaluated at interpolation nodes.

Before describing the implementation details for the Fourier and RBF models, we discuss the node and evaluation points used.

\subsection{Node and Evaluation points}
For our 2D objects, we use $N$ equally-spaced points on the interval
$(-\pi,\pi]$ as the node set $\{\lon_k\}_{k=1}^{N}$, and take $N$ to
be even.  This gives a uniform sampling in the parametric space and
allows fast algorithms to be used for computing the interpolants as
discussed below.  Additionally, since the shape of our target objects are near circular or elliptical, these nodes give a good distribution of data sites on the object.  We also use $M >> N$
equally-spaced points in the interval $(-\pi,\pi]$ as the set of evaluation
points $\{\lons_j\}_{j=1}^{M}$ since this also results in a set of sample sites that are well distributed over the object.

To get a good sampling of our nearly ellipsoidal or
spherical objects in 3D, we cannot resort to using equally spaced
points in the spherical coordinate system as our node sets
$\{(\lon_k,\lat_k)\}_{k=1}^{N}$ because of the inherent ``pole
problem''.  Instead we use node sets that give a quasi-uniform
distribution of data sites on the unit sphere.  Since only a maximum
of 20 points can be evenly distributed on a sphere, there are a myriad
of methods to define and generate a quasi-uniform distribution for
larger numbers of points~\cite{HardinSaff:2004}.  We use two of these
methods: maximal determinant (MD) for our spherical harmonic models
and minimal energy (ME) for our RBF models. Both of these methods are
discussed in~\cite{WomersleySloan:2001} and many of these
two point sets for various $N$ can be downloaded
from~\cite{WomerSloan:2003}.  The MD points are generated by finding a
distribution of points that maximize the determinant of a certain
``Gram matrix'' related to \eqref{eq:sph_linsys}.  The ME points are
generated by finding a distribution of nodes that minimize an
electrostatic type energy potential.  For spherical harmonic
interpolation, the MD points lead to much better results both in terms
of accuracy and stability~\cite{WomersleySloan:2001}.  For RBF
interpolation, the ME points typically yield better results in terms
of accuracy~\cite{FlyerWright:2007,FlyerWright:2009} for larger shape
parameters $\ep$.  For smaller values, the MD points give better
results because of the connection to spherical harmonics as
$\ep\rightarrow 0$~\cite{FornbergPiret:2007}.  For the set of
evaluation points, $\{(\lons_j,\lats_j)\}_{j=1}^M$, we use $M >> N$ ME
points for both the spherical harmonic and RBF models, which again
results in a well distributed set of sample sites on the object.

\subsubsection{Fourier Models}\label{sec:implementation_trig}
The first step in computing the normal vectors and forces for the 2D
trigonometric model \eqref{eq:2D_obj_trig} is to compute the
interpolation coefficients $c_k^x$ and $c_k^y$, $k=1,\ldots,N$ (see
\eqref{eq:trig_interp}).  Since we are using equally spaced node
points $\{\lond_k\}_{k=1}^N$, we can avoid having to solve
\eqref{eq:trig_linsys} directly for these coefficients and can instead
compute them by means of the fast Fourier transforms ({\em
  e.g.}~\cite[\S 3]{Trefethen:2000:SMM}) at a cost of $O(N\log N)$.

We next compute the derivatives of the interpolants to obtain the
following approximation to \eqref{eq:tngt_2d}: 
\begin{align}
\frac{\partial}{\partial \lon}\vx(\lon)\bigr|_{\lon = \lons_j} \approx \frac{\partial}{\partial \lon}{\vp}(\lon)\bigr|_{\lon = \lons_j},\;j=1,\ldots,M.
\label{eq:trig_deriv1}
\end{align}
We then determine the normal vector at $\vx(\lons_j)$ by normalizing
the vector above and switching the components according to
\eqref{eq:utngt_2d} and \eqref{eq:unnrml_2d}. We similarly obtain an
approximation of the force \eqref{eq:force_2d} from the second
derivative of the interpolants:

\begin{align}
  \frac{\partial^2}{\partial \lon^2}\vx(\lon)\bigr|_{\lon = \lons_j}
  \approx \frac{\partial^2}{\partial
    \lon^2}{\vp}(\lon)\bigr|_{\lon =
    \lons_j},\;j=1,\ldots,M.
\label{eq:trig_deriv2}
\end{align}

For the 3D spherical harmonic model \eqref{eq:3D_obj_sph}, the first
step in computing the normal vectors and forces is again to compute
the interpolation coefficients $c_k^x$, $c_k^y$, and $c_k^z$,
$k=1,\ldots,N$, (see \eqref{eq:sph_interp}).  Unlike the trigonometric
interpolant, there are unfortunately no fast algorithms for computing
these coefficients.  Since we use relatively small values of $N$, we
thus resort to determining the coefficients by solving the linear
system \eqref{eq:sph_linsys} using a direct $LU$ factorization of the
interpolation matrix.  By using the MD points as the nodes in this
model, we are guaranteed that this system is non-singular and
relatively well conditioned~\cite{WomersleySloan:2001}.  We note that
in context of the IB method simulation, the node points will stay
fixed throughout the simulation so that the $LU$ factorization of the
interpolation matrix from \eqref{eq:sph_linsys} needs to be done
only once at the initial time-step and then stored.  Thus,
for all other time-steps the coefficients can be determined in
$O(N^2)$ computations.

After the coefficients are determined, we compute the following six
derivatives to obtain approximations to \eqref{eq:tngt1_3D} and
\eqref{eq:tngt2_3D}:
\begin{align}
\frac{\partial}{\partial \lon}\vx(\lon,\lat)\bigr|_{(\lon,\lat) = (\lons_j,\lats_j)} & \approx \frac{\partial}{\partial \lon}{\vp}(\lon,\lat)\bigr|_{(\lon,\lat) = (\lons_j,\lats_j)},\;j=1,\ldots,M, 
\label{eq:sph_deriv1}\\
\frac{\partial}{\partial \lat}\vx(\lon,\lat)\bigr|_{(\lon,\lat) = (\lons_j,\lats_j)} & \approx \frac{\partial}{\partial \lat}{\vp}(\lon,\lat)\bigr|_{(\lon,\lat) = (\lons_j,\lats_j)},\;j=1,\ldots,M. 
\label{eq:sph_deriv2}
\end{align}
We then compute the normal vectors using these approximations in \eqref{eq:utngt_3d} and \eqref{eq:unrml_3d}. 

The computation of the force requires the approximation to the normal
vectors and an approximation to the mean curvature
\eqref{eq:mean_curvature}.  For the values of $E$, $F$, and $G$ in the
mean curvature computation (see \eqref{eq:first_fund_form}), we use the
approximations \eqref{eq:sph_deriv1} and \eqref{eq:sph_deriv2}.  For
the values of $e$, $f$, and $g$, we use the approximations
\begin{align}
\frac{\partial^2}{\partial \lon^2}\vx(\lon,\lat)\bigr|_{(\lon,\lat) = (\lons_j,\lats_j)} & \approx \frac{\partial^2}{\partial \lon^2}{\vp}(\lon,\lat)\bigr|_{(\lon,\lat) = (\lons_j,\lats_j)},\;j=1,\ldots,M, 
\label{eq:sph2_deriv1}\\
\frac{\partial^2}{\partial \lat \partial \lon}\vx(\lon,\lat)\bigr|_{(\lon,\lat) = (\lons_j,\lats_j)} & \approx \frac{\partial^2}{\partial \lat \partial \lon}{\vp}(\lon,\lat)\bigr|_{(\lon,\lat) = (\lons_j,\lats_j)},\;j=1,\ldots,M, 
\label{eq:sph3_deriv1}\\
\frac{\partial^2}{\partial \lat^2}\vx(\lon,\lat)\bigr|_{(\lon,\lat) = (\lons_j,\lats_j)} & \approx \frac{\partial^2}{\partial \lat^2}{\vp}(\lon,\lat)\bigr|_{(\lon,\lat) = (\lons_j,\lats_j)},\;j=1,\ldots,M. 
\label{eq:sph2_deriv2}
\end{align}

\subsubsection{RBF Models}\label{sec:implementation_rbf}
The normal vectors and forces for the RBF models are computed in the
same fashion as for the Fourier models discussed above; one
just needs to replace the Fourier interpolants
${\vp}(\lon)$ and ${\vp}(\lon,\lat)$ with
the RBF interpolants ${\vs}(\lon)$ from
\eqref{eq:2D_obj_rbf} and ${\vs}(\lon,\lat)$ from
\eqref{eq:3D_obj_rbf}, respectively.  We thus omit a full description.
We will, however, discuss the shape parameter $\ep$ and the
computation of the interpolation coefficients.

Infinitely smooth radial kernels like the MQ \eqref{eq:mq} and IMQ
\eqref{eq:imq} feature a free shape parameter $\ep$.  It has generally
been reported in the literature that there is typically an optimal
value of $\ep$ that produces the best accuracy in the interpolants
with these kernels and that this value tends to decrease with
increasing smoothness of the underlying function being approximated
({\em e.g.}~\cite{Rippa:1999}).  However, as $\ep$ decreases to zero
these smooth kernels become increasingly flat and the shifts of $\phi$
in \eqref{eq:circ_rbf_interp} and \eqref{eq:sph_rbf_interp} become
less and less distinguishable from one another.  If one follows the
direct approach of solving for the expansion coefficients via
\eqref{eq:rbf_linsys} and then evaluating the interpolant via
\eqref{eq:circ_rbf_interp} or \eqref{eq:sph_rbf_interp} (which is
denoted by RBF-Direct in the current literature) for $\ep$ in this
flat regime, then ill-conditioning can entirely contaminate the
computation.  For RBF interpolation on a sphere, this ill-conditioning
can be completely bypassed by replacing the RBF-Direct algorithm with
the RBF-QR algorithm of Fornberg and Piret~\cite{FornbergPiret:2007}.
The framework for this algorithm can also naturally be adapted to the
task of computing RBF interpolants on the unit circle in a stable
manner for all $\ep$.

We have implemented both the RBF-QR algorithm and the RBF-Direct
approach and present results in Sections \ref{sec:results2d_shape} and
\ref{sec:results3d_shape} illustrating the behavior of the RBF
interpolants for the full range of $\ep$ and the connection to Fourier
based methods.  However, we have opted to use the RBF-Direct approach
in implementation since it is computationally more efficient and the
coding is much less involved for computing the normals and forces.
Additionally, we have found that with the RBF-Direct approach and the
values of $N$ that we considered it is possible to get as good or
better results than the Fourier based methods.  For increasingly large
values of $N$, or objects whose parameterizations are very smooth, it
may be necessary to switch to the RBF-QR algorithms to exploit the
better accuracy that can sometimes be achieved for increasingly small
values of $\ep$.

For the RBF-Direct approach, the interpolation coefficients for both
the 2D and 3D objects can be determined by solving the
linear system \eqref{eq:rbf_linsys} (with the appropriate choice of
$r_{j,k}$ for the dimension of interest).  In the case of 2D objects with equally spaced points,
solving this system directly can be bypassed by means of the fast
Fourier transform and the coefficients can be computed in $O(N\log N)$
operations~\cite{HubbertMuller:2006}.  This follows by observing that
the matrix in \eqref{eq:rbf_linsys} is \emph{circulant} (for any
radial kernel $\phi$) and can be diagonalized via the discrete Fourier
transform matrix~\cite[\S 4.7.7]{GolubVanLoan:1996}.  For the 2D
models, we use the MQ radial kernel \eqref{eq:mq}.

As in the case of the spherical harmonic model, there are no fast
direct algorithms for determining the interpolation coefficients for the 3D
RBF model \eqref{eq:3D_obj_rbf} and we thus resort to using a direct
method.  However, unlike the spherical harmonic model the system is
symmetric and, as discussed in Section \ref{sec:modeling_rbf}, for the
right choice of $\phi$ it is positive definite.  Thus, a Cholesky
factorization of the matrix can be used which reduces the memory costs
and the need for pivoting over the $LU$ factorization method used in
the spherical harmonic model.  We also note that the initial cost of
computing the Cholesky factorization is lower than the $LU$
factorization, but since this is only done once initially there is no
real savings in an IB simulation.  We have opted to use the IMQ
kernel \eqref{eq:mq} to exploit the use of the Cholesky factorization.

\section{2D Platelet Modeling Results}
\label{sec:Results2D}

In this section, we present the results of our comparative study
between using the piecewise linear approach as traditionally
used within the IB method and our two alternative parametric
approaches in 2D: RBF and Fourier (trigonometric polynomials)
interpolation.  Recall that within an IB timestep, the typical
procedure employed is as follows.  Given the locations of the immersed
boundaries, both the normals and forces on an object are computed. The
forces are then projected to an Eulerian grid and used as
right-hand-side forcing to the Navier-Stokes equations.  Based upon an
update velocity field, the positions of the IB points are updated.  In
our comparison, we thus examine the geometric modeling capabilities,
accuracy of the normal computations, and accuracy of the computation
of the forces.  As discussed in the previous section, we distinguish
between the data sites and sample sites for the parametric models.
Data sites are the positions along the object at which the parametric
models are interpolating. It is at these positions that we propose
updating the geometric information of the object (for instance, at the
conclusion of a timestep when the object's movement within
the flow field is updated).  Sample sites (which are normally more
numerous compared to the data sites) are the positions along the
object at which normals and forces are computed.  It is from these
positions that we propose projecting the IB forces.  In all
experiments, $100$ sample sites are used as this represents the
typical number of IB points that would be used per platelet object in
a traditional 2D immersed boundary computation (and hence a reasonable
standard against which to compare our new methods for the purposes of
determining the feasibility of replacement).  All errors are computed
by taking the maximum of the two-norm difference between the
approximations and the true values.

\subsection{Test Cases}
\label{sec:results2d_testcases}

We consider 2 prototypical test objects and define them based upon
perturbations of idealized shapes (an ellipse and a circle).  Let
\(\textbf{x}_{ideal}\) be a function representing the idealized,
unperturbed shapes as given by the following equation:
\begin{equation}
\vx_{ideal} = (x_c + a \cos \lon, y_c + b \sin \lon)
\end{equation}
where $-\pi \leq \lon \leq \pi$.  Here $(x_c,y_c)$ denotes the object
center and $a$ and $b$ denote the radii.  The two objects used for our
comparison are defined as follows:
\begin{eqnarray}
\text{\underline{Object 1}:}\quad \vx_{2d\,obj 1} &=& \lf[1.0 + A \exp\lf(\frac{-(1-\cos \lon)^2}{\sigma_1}\rt)\rt] \, \vx_{ideal}, \label{eq:obj1_2d}\\
\text{\underline{Object 2}:}\quad \vx_{2d\,obj 2} &=& \lf[1.0 + B \exp\lf(\frac{(-(1-\cos^{2}\lon)^{1.5})}{\sigma_2}\rt)\rt] \, \vx_{ideal} \label{eq:obj2_2d}.
\end{eqnarray}
For Object 1, we use the following parameters: $x_c = y_c = 0.9$,
$a=0.04$, $b=0.05$, $A = 0.09$ and $\sigma_1 = 0.1$.  For Object 2, we
use the following parameters: $x_c = y_c = 0.2$, $a=b=0.1$, $B = 0.04$
and $\sigma_2 = 0.9$.
 
\begin{figure}[htbp]
\begin{center}
\includegraphics[width=6.0in]{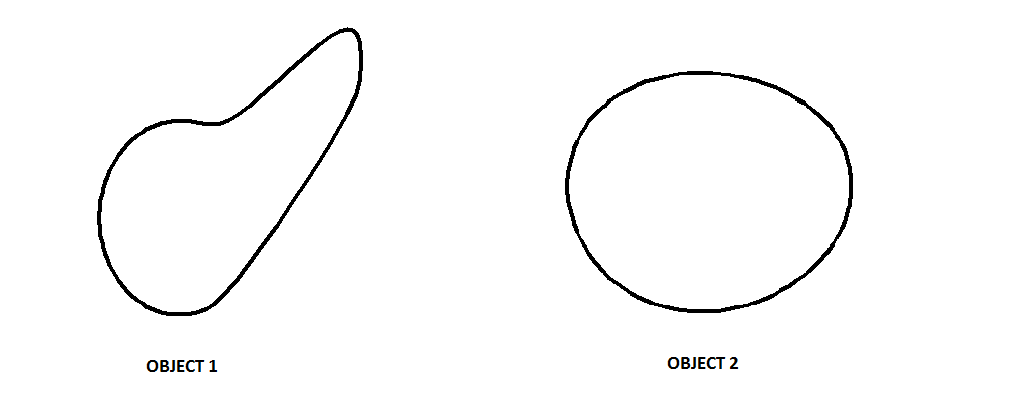}
\end{center}
\caption{The test objects \eqref{eq:obj1_2d} and \eqref{eq:obj2_2d} for the 2D study.} 
\label{fig:diagram2d}
\end{figure}

Figure \ref{fig:diagram2d} displays the two test objects
\eqref{eq:obj1_2d} and \eqref{eq:obj2_2d}.  Object 1 is a smooth (in
terms of regularity) yet highly perturbed ellipse, while Object 2 is a
non-smooth perturbation of a circle.  It can be shown that the
parameterization \eqref{eq:obj2_2d} for this object has only
two continuous derivatives.

\subsection{Comparison of Reconstructing the Objects}
\label{sec:results2d_geoerror}
\begin{figure}[htbp]
\begin{center}
\includegraphics[width=3.0in]{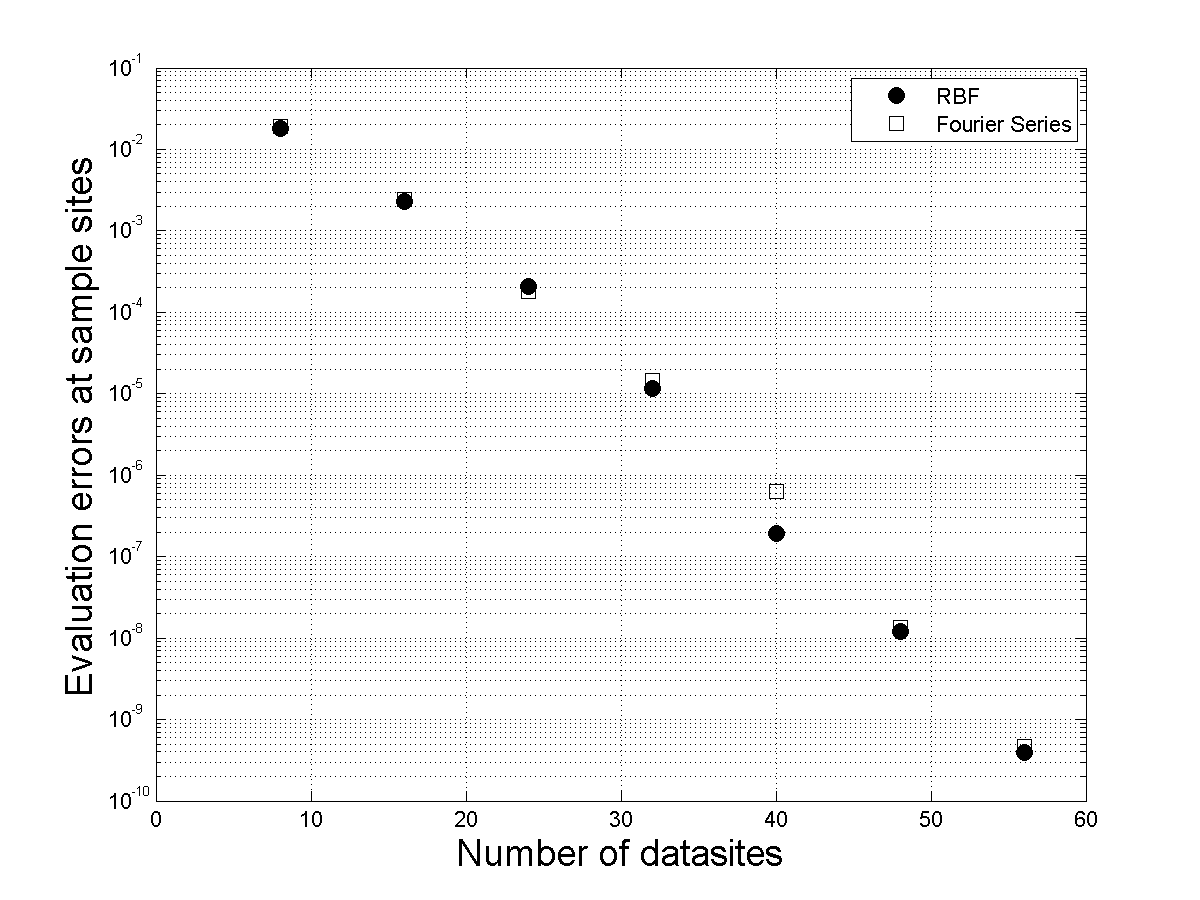}
\includegraphics[width=3.0in]{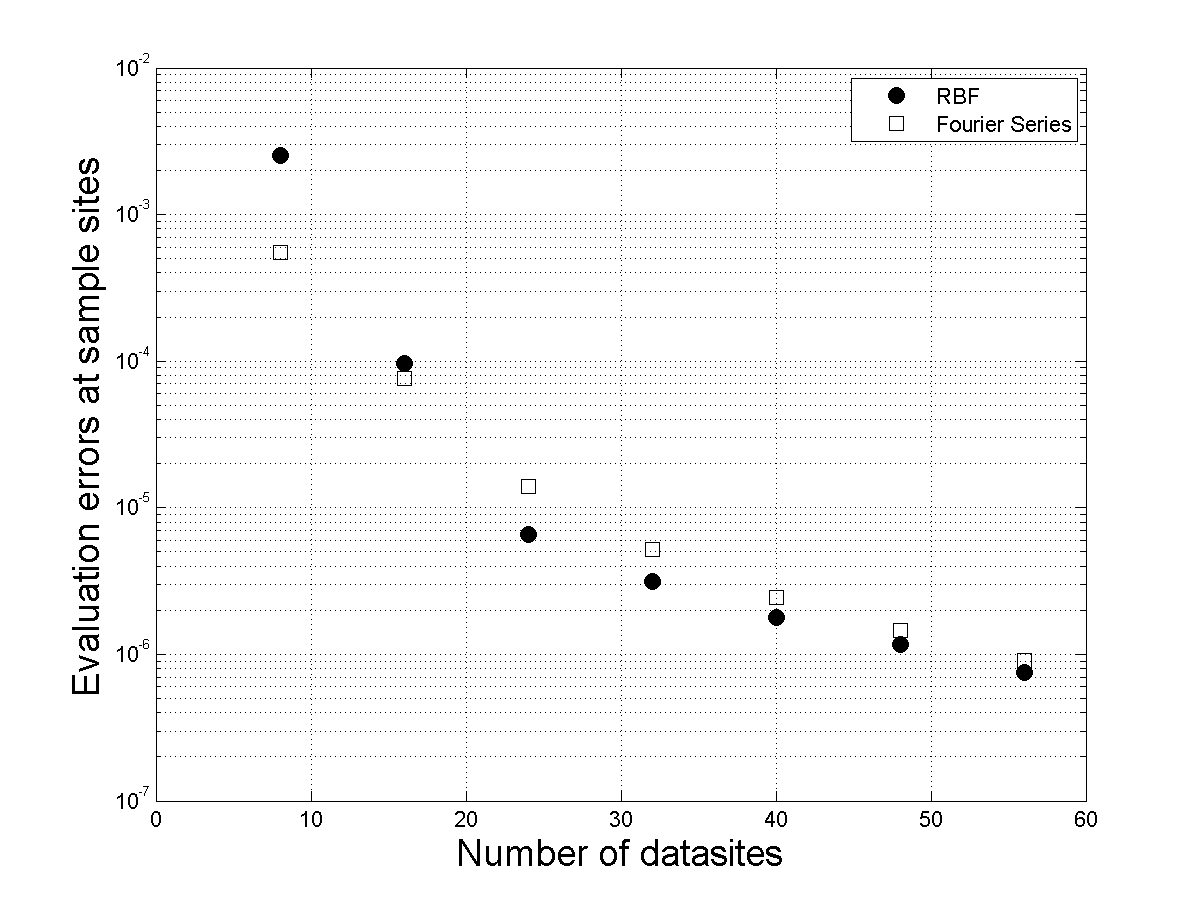}
\end{center}
\caption{Error in the reconstruction of the shape of the objects (left
  is Object 1 and right is Object 2) evaluated at $M=100$ sample sites
  as a function of the number of data sites. Circles denote the errors
  in the RBF model and squares denote the errors for the Fourier
  model. For the RBF model, the shape parameter for Object 1 was set
  to $\ep=0.9$ and for Object 2, it was set to $\ep=3.6$.}
\label{fig:geom2d}
\end{figure}

We first examine the errors in reconstructing the objects using the
RBF and Fourier approaches. In Figure \ref{fig:geom2d} we present the
errors in reconstructing the objects as a function of the number of
data sites.  The error at the sample sites gives an indication of the
modeling capability of the RBF and Fourier methods.  We can see from
this figure that both the RBF and Fourier models are converging at a
spectral rate for Object 1 (left figure), but at a much slower rate
for Object 2 (right figure).  This is expected since Object 1 is
infinitely smooth, while Object 2 has only two continuous derivatives.
The RBF and Fourier models perform similarly for Object 1.  For Object
2, the RBF model shows better reconstruction properties as the number
of sample sites increases above 20.  No direct comparison with the
piecewise linear model is given as the piecewise linear IB method
always samples at the interpolating points.  

\subsubsection{Shape Parameter Study}
\label{sec:results2d_shape}
\begin{figure}[htbp]
\begin{center}
\includegraphics[width=3.0in]{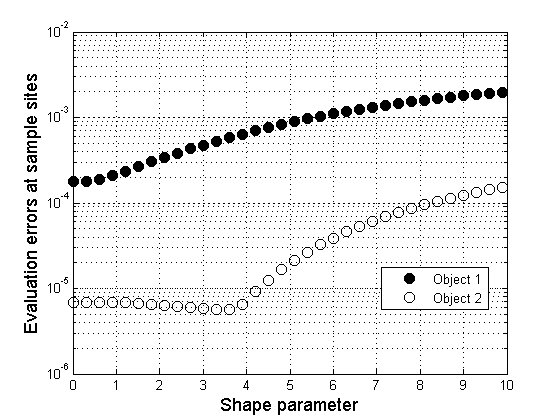}
\end{center}
\caption{Errors in the RBF reconstructions of the objects using $N=24$ data sites and $M=100$ sample sites as a function of the shape parameter.} 
\label{fig:2dshape1}
\end{figure}

\begin{figure}[htbp]
\begin{center}
\includegraphics[width=3.0in]{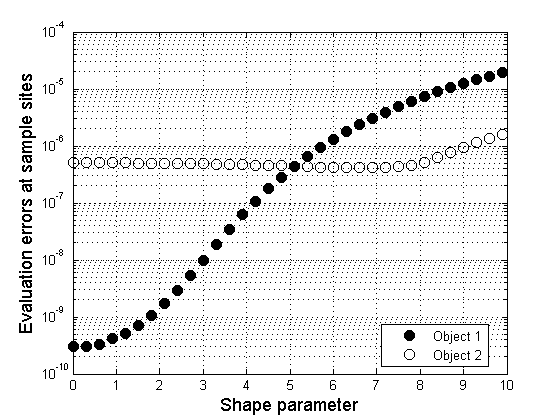}
\end{center}
\caption{Errors in the RBF reconstructions of the objects using $N=56$
  data sites and $M=100$ sample sites as a function of the shape
  parameter.}
\label{fig:2dshape2}
\end{figure}

In this section, we examine the impact of the shape parameter on the
reconstruction errors of the RBF model.  Figure \ref{fig:2dshape1}
displays the reconstruction errors for the two objects as a function
of the shape parameter using $N=24$ data sites.  A similar comparison
for $N=56$ data sites is given in Figure \ref{fig:2dshape2}.  For $\ep
\lesssim 0.85$, it was necessary to use the RBF-QR
algorithm~\cite{FornbergPiret:2007} (adapted to the unit circle) to
compute the model in a numerically stable manner for the $N=56$ case.
We can see from both figures that the errors are smallest for
$\ep\approx 0$ for the smooth Object 1 and increase quite dramatically
as $\ep$ increases.  For the non-smooth Object 2, there is a much
larger range of $\ep$ for which the errors are small, and this range
includes values for which the RBF-Direct approach can be used without
issues of numerical instabilities.

We used Figures \ref{fig:2dshape1} and \ref{fig:2dshape2} to
help guide our selection of $\ep$ for the numerical experiments.
However, we found from extensive tests on other objects, that if the
object is smooth, and RBF-Direct is to be used, then one generally
wants to choose $\ep$ as small as the numerical conditioning allows.
For non-smooth objects there is much more freedom in the choice and
the results will not vary that greatly.  It is unclear if we should
expect smooth or non-smooth objects in an IB simulation.

We conclude this section by noting that there are several algorithms
that have been devoted to selecting an ``optimal'' shape
parameter~\cite[\S 17]{Fasshauer:2007}.  However, these are too costly
to be used every time-step of an IB simulation.  We are thus
advocating using a fixed $\ep$ for all time-steps.  This value could
be selected based on an expected typical shape for the immersed
objects and one of the algorithms from~\cite[\S 17]{Fasshauer:2007} or
from trial and error.  We will report on these strategies in a follow
up paper where the RBF models are used in actual IB method
simulations.

\subsection{Comparison of Normal Vectors and Forces}
\label{sec:results2d_force}
We next focus on the errors in the parametric models in the
approximation of the normal vectors to the objects and the forces.  In
this case, we compare the results to the traditional piecewise linear
models.

\begin{figure}[htbp]
\begin{center}
\includegraphics[width=3.0in]{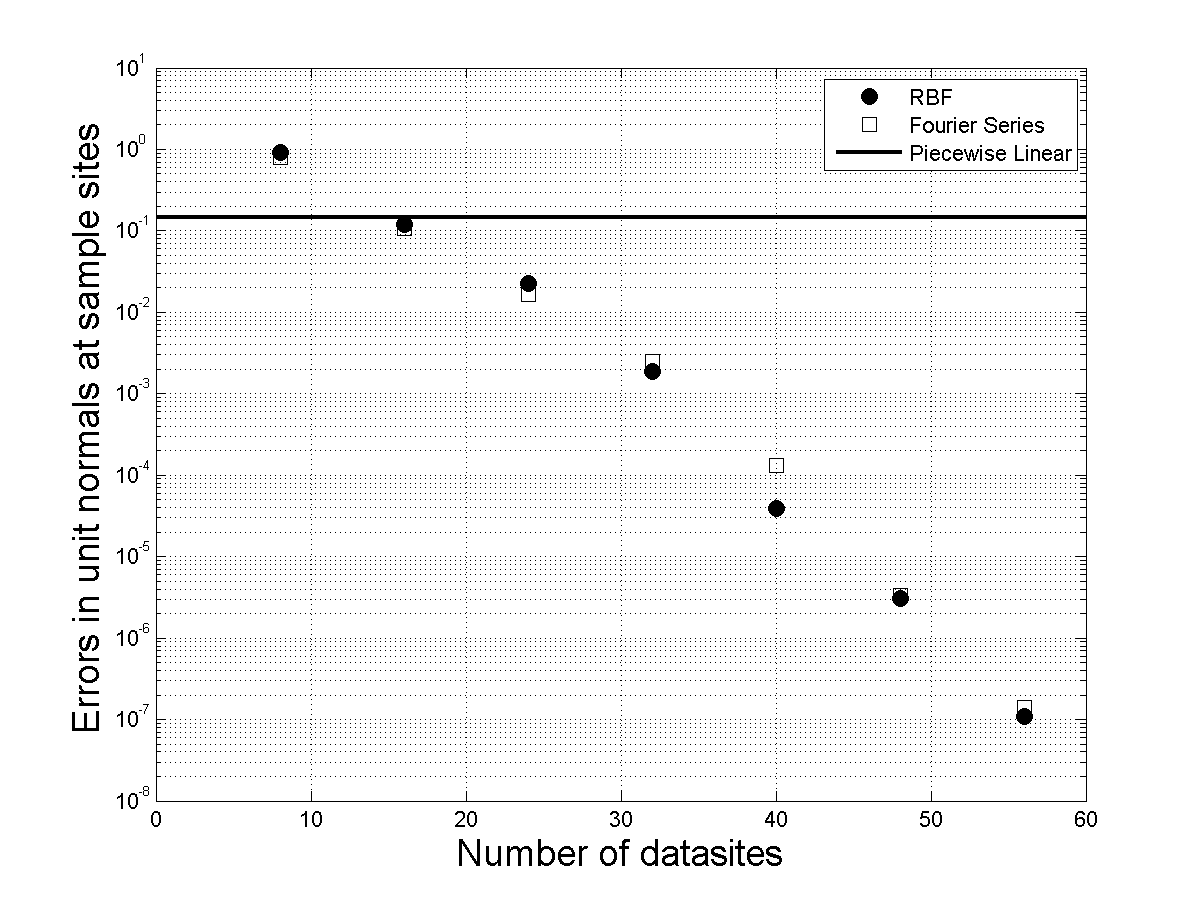}
\includegraphics[width=3.0in]{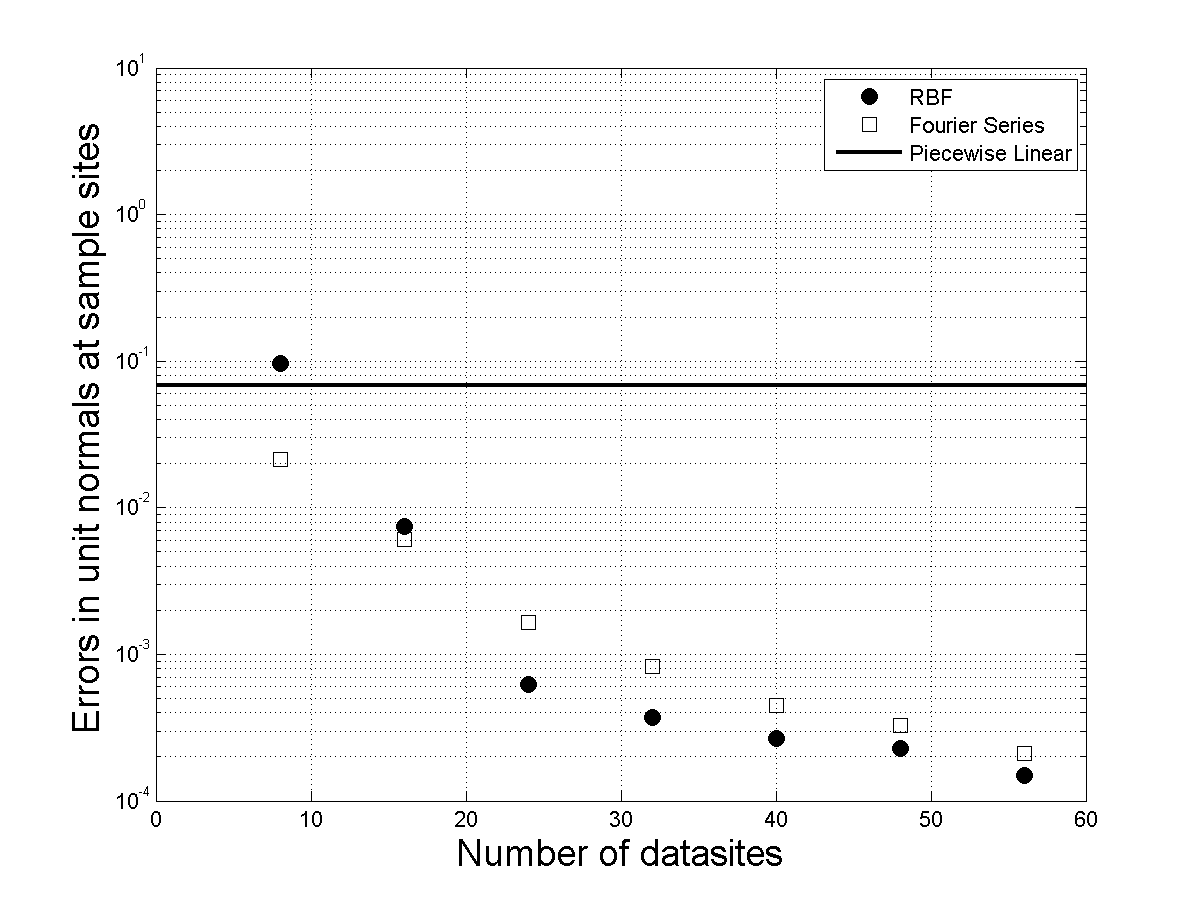}
\end{center}
\caption{Errors in the approximations of the normal vectors to the
  objects at $100$ sample sites as a function of the number of
  data sites $N$.  The left plot is for Object 1, while the right one
  is for Object 2. The line denotes the error for the method used in
  the piecewise linear model with $100$ IB points. Circles denote the
  errors for the RBF model and squares denote the Fourier model. For
  the RBF model, $\ep=0.9$ for Object 1 and $\ep=3.6$ for Object 2.}
\label{fig:normals2d}
\end{figure}

Figure \ref{fig:normals2d} displays the errors in the normal vectors
at $M=100$ sample sites as a function of the number of data sites $N$
for both the RBF and Fourier models.  A solid line denoting the errors
in the normal vectors for $100$ IB points is given for comparison
using the method for the piecewise linear models discussed in Section
\ref{sec:implementation_pwl}. We can see from this figure that at
about $N=18$ data sites the errors for both the RBF and Fourier models
of Object 1 are lower than the piecewise linear model.  The errors are
similar between both parametric models and decrease rapidly with
increasing $N$.  The results for Object 2 are even more favorable for
the parametric models compared to the piecewise linear model.  For
increasing $N$ the RBF model appears to have an advantage over the
Fourier model.

\begin{figure}[htbp]
\begin{center}
\includegraphics[width=3.0in]{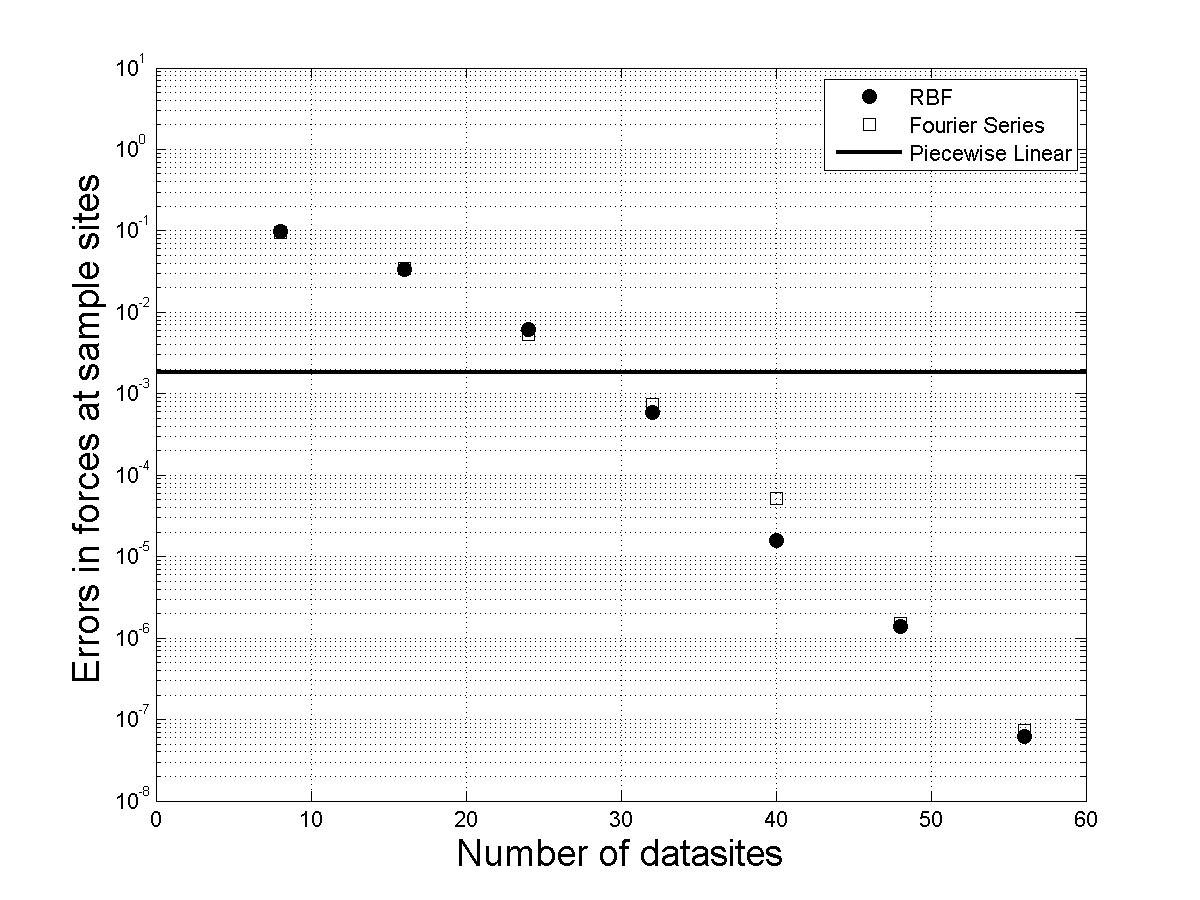}
\includegraphics[width=3.0in]{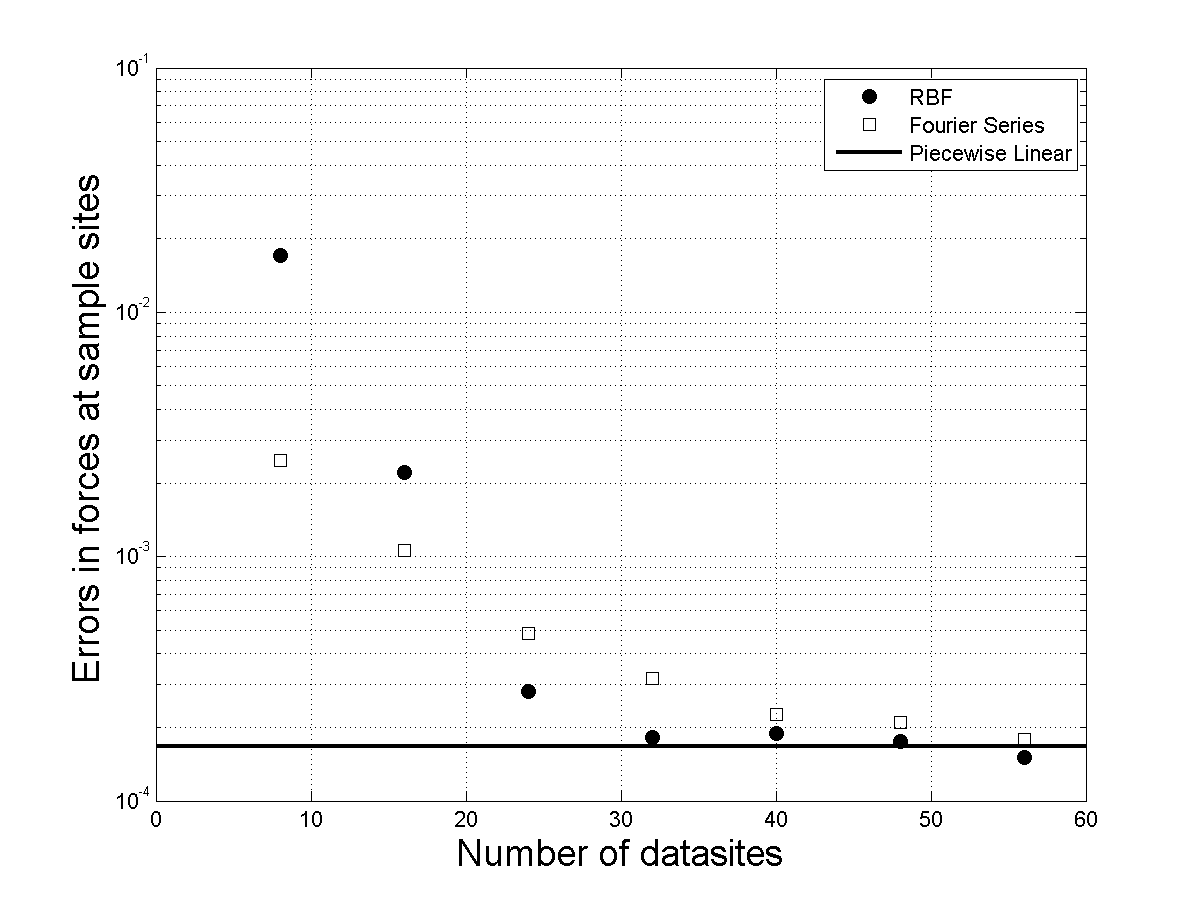}
\end{center}
\caption{Errors in the approximation of the forces evaluated at
  $M=100$ sample sites as a function of the number of data sites $N$
  for Object 1 (left) and Object 2 (right). The black line denotes the
  errors for a piecewise linear model with $100$ IB points. Circles
  denote the errors for the RBF model and squares denote the errors
  for the Fourier model. For the RBF model, $\ep=0.9$ for Object 1 and $\ep=3.6$ for
  Object 2.}
\label{fig:force2d}
\end{figure}

We lastly examine the errors in the force computation incurred by the
two parametric models and the traditional piecewise linear
model. Figure \ref{fig:force2d} shows the errors in forces evaluated
at $100$ sample sites as a function of the number of data sites
$N$. In all experiments, the force constant $K_0$ is set to $0.2$.
The solid line in Figure \ref{fig:force2d} denotes the error for the
piecewise linear model computed at $100$ IB points.  For Object 1, we
can see from the left plot of this figure that the errors for both
parametric models are lower than the piecewise linear model starting
at about $N=30$ data sites.  Again, both the RBF and Fourier models
give similar results for this object.  For the non-smooth Object 2, it
requires about $N=32$ data sites for the RBF model to match the errors
of the piecewise linear model, while it takes approximately $N=56$
data sites for the Fourier model to give similar errors. We note
that the errors for both the RBF and Fourier models do not fall as sharply
for the non-smooth Object 2 as the number of datasites is increased. This is
because Object 2 is generated from a function that has only two derivatives,
and the force computation involves computing a second derivative. It therefore
follows that these global methods would therefore not converge as they would
in the case of the smooth Object 1.

\subsection{Comparison of the Computational Cost}
\label{sec:results2d_compcost}

We conclude the 2D results experiments with an examination of the
computational cost associated with the three methods.  We measure the
computational cost as the elapsed wallclock time required to compute
the interpolation coefficients, evaluate the interpolants, compute the
normal vectors and compute the forces. Under the assumption that all
objects will be evaluated at the same parametric sites at each
timestep, for both parametric models we pre-compute the matrices for
evaluating the interpolants, the derivatives, and the force operator
once the interpolation coefficients have been determined (see Section
\ref{sec:implementation_param} for details).  We do not account for
this setup time in our timing results.

Since for the piecewise linear model the number of evaluation sites is
the same as the number of data sites, the total computational cost
includes only the time required to compute the normal vectors and
forces (see Section \ref{sec:implementation_pwl} for details).

All computations were performed in \matlab version 7.10.0499 (64-bit)
on a Windows desktop with a Intel Core i7 Sandy Bridge 3.4 GHz
processor and 4 GB of 1600 MHz RAM.  Times were measured using the tic
and toc functions in \matlab. All results presented are averages of a
100 trials and are in seconds.

\begin{figure}[htbp]
\begin{center}
\includegraphics[width=3.0in]{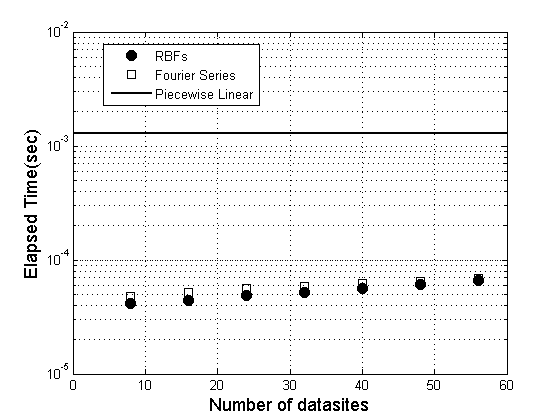}
\end{center}
\caption{Elapsed wallclock time (in seconds) for one object to perform
  interpolation, evaluation, the computation of normal vectors, and
  the computation of forces at $M=100$ sample sites as a function of
  the number of data sites $N$.  The piecewise linear computations
  were done with $100$ IB points for comparison.}
\label{fig:results2d_compcost}
\end{figure}

Figure \ref{fig:results2d_compcost} displays the elapsed time between
the RBF, Fourier, and traditional piecewise linear models. The results
for the RBF and Fourier models are displayed as a function of the
number of data sites $N$ for a fixed number of $M=100$ sample sites.
The results for the piecewise linear model are for a fixed number of
$100$ IB points.  We can see from the figure that the parametric
models require significantly less time than the piecewise linear
model. For $N=56$ data sites, the parametric models are over one order
of magnitude faster.

We note that all the evaluation and derivative computations for the
parametric models can be formulated in terms of matrices in order to
avoid the need to first solve for the coefficients every time step of
the IB simulation.  Thus, the results we present are not optimal in
terms of computational time.  If, however, during the IB simulation
the sample sites change, then the step of going first through the
coefficients as we have done will be necessary.

\section{3D Platelet Modeling Results}
\label{sec:Results3D}

Following a similar approach to the last section, we present here the results from a comparative study between using the traditional piecewise linear approach as used within the IB method and our two alternative parametric approaches in 3D: RBF and Fourier (spherical harmonics) interpolation.  We examine the reconstruction capabilities of the models and the accuracy in computing normal vectors and forces.  As in the 2D tests, we distinguish between data sites and sample sites. In all experiments unless otherwise specified, $M=1024$ sample sites are used as this represents the typical number of IB points that would be used per platelet object in a traditional 3D IB computation (and hence a reasonable standard against which to compare our new methods for the purposes of determining the feasibility of replacement). All errors are computed by taking the maximum of the two-norm difference between the approximations and the true values.

\subsection{Test Cases}
\label{sec:results3d_testcases}

We again consider 2 prototypical test objects and define them based on perturbations of idealized shapes (an ellipsoid and a sphere).  Let \(\textbf{x}_{ideal}\) be a function representing the idealized, unperturbed shapes as given by the following equation:
\begin{equation}
\vx_{ideal} = (x_c + a \cos \lon \cos \lat, y_c + b \sin \lon \cos \lat,  z_c + c \sin \lat),
\end{equation}
where $-\pi \leq \lon \geq \pi$ and $-\frac{\pi}{2}\leq \lat \leq \frac{\pi}{2}$.  Here $(x_c,y_c,z_c)$ denotes the object center, $a$ and $b$ are the equatorial radii, and $c$ is the polar radius. The two objects used for our comparison are defined as follows:
\begin{eqnarray}
\text{\underline{Object 1}:  }\quad \vx_{3d\,obj 1} &=& \lf[1.0 + A \exp\lf(\frac{r_c^2}{\sigma_1}\rt)\rt]\,\vx_{ideal}, \label{eq:obj1_3d}\\
\text{\underline{Object 2}:  }\quad \vx_{3d\,obj 2} &=& \lf[1.0 + B \exp\lf(\frac{r_c^{2.5}}{\sigma_2}\rt)\rt]\,\vx_{ideal}. \label{eq:obj2_3d}
\end{eqnarray}
where $r_c = 1-\cos \lat \cos \lat_c \cos(\lon-\lon_c)-\sin \lat \sin \lat_c$. 
For Object 1, we use the following parameters: $x_c = y_c = z_c = 0.9$, 
$a=0.1$, $b=0.2$, $c=0.09$, $A = 0.09$ and $\sigma_1 = 0.2$.  
For Object 2, we use the following parameters: 
$x_c = y_c = 0.1$, $z_c = 0.2$, $a = b= c = 0.1$, $B = 0.04$ and $\sigma_2 = \frac{16}{25}$.
For both objects $\lon_c = 0$ and $\lat_c = \frac{\pi}{2}$.

\begin{figure}[htbp]
\begin{center}
\includegraphics[width=6.0in]{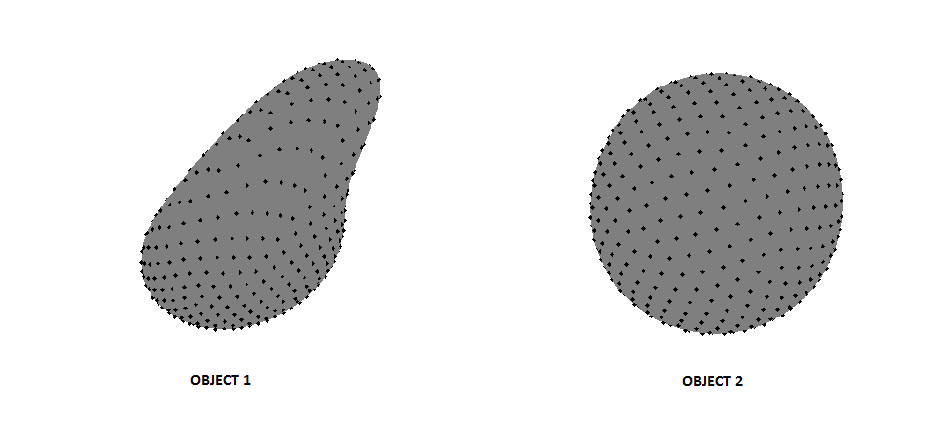}
\end{center}
\caption{The test objects \eqref{eq:obj1_3d} and \eqref{eq:obj2_3d} for the 3D study.} 
\label{fig:diagram3d}
\end{figure}

Figure \ref{fig:diagram3d} displays the two test objects
\eqref{eq:obj1_3d} and \eqref{eq:obj2_3d}.  Object 1 is a smooth (in
terms of regularity) yet highly perturbed ellipsoid, while the Object
2 is a non-smooth perturbation of a sphere.  It can be shown that the
parameterization \eqref{eq:obj2_3d} has only three continuous derivatives.

\subsection{Comparison of Reconstructing the Objects}
\label{sec:results3d_geoerror}
\begin{figure}[htbp]
\begin{center}
\includegraphics[width=3.0in]{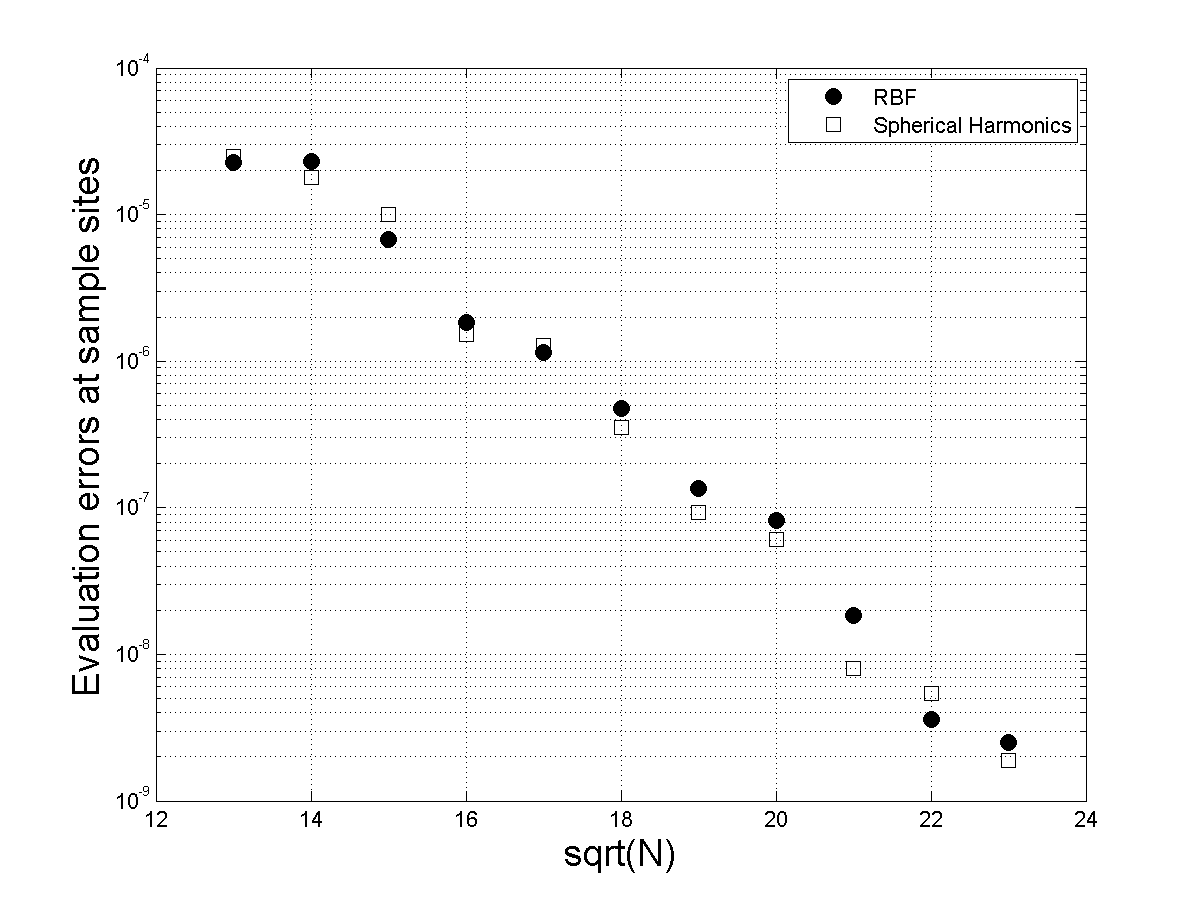}
\includegraphics[width=3.0in]{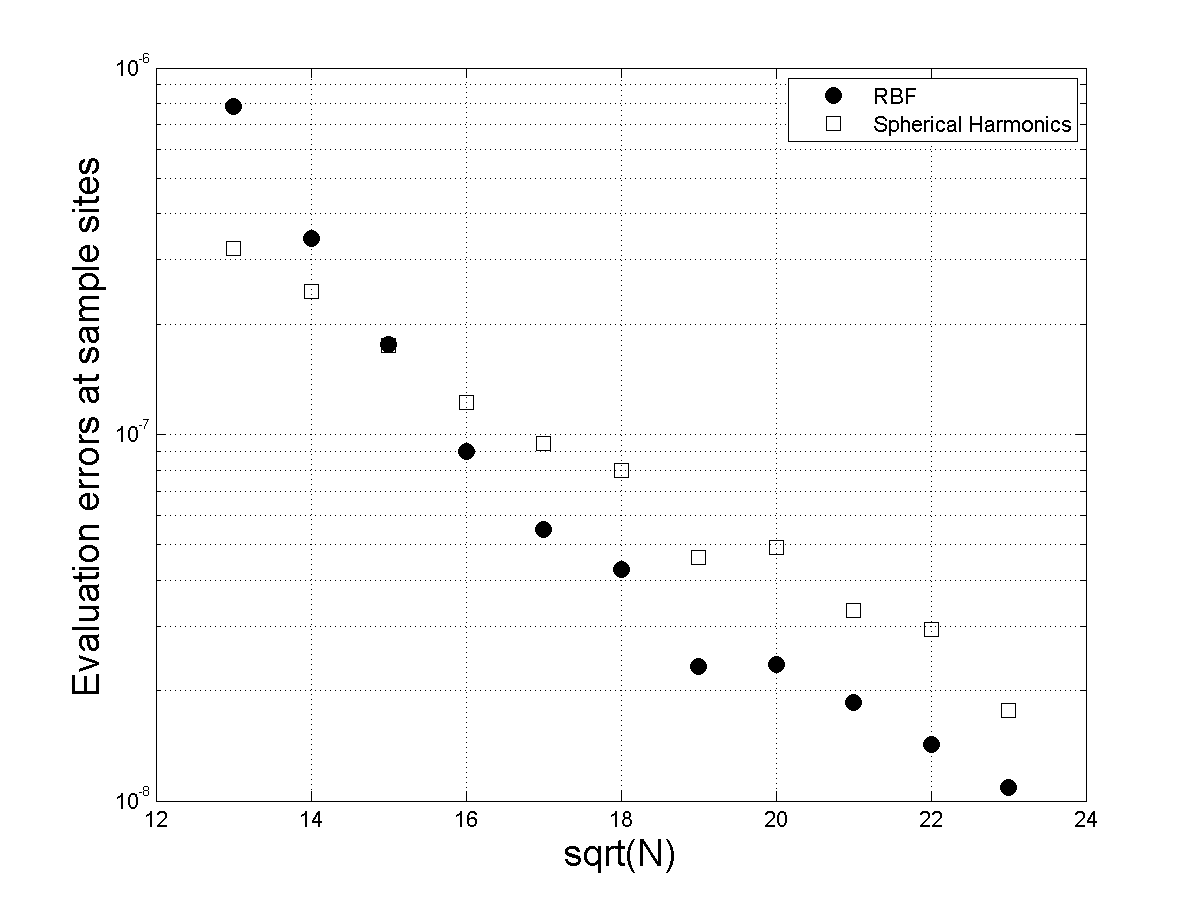} 
\end{center}
\caption{Error in the reconstruction of the shape of the objects (left is Object 1 and right is Object 2) evaluated at $M=1024$ sample sites
as a function of the square root of the number of data sites $N$. Circles denote the errors in the RBF model and squares denote the errors for the Fourier model. For the RBF model, the shape parameter for Object 1 was set to $\ep=0.9$ and for Object 2 was 
set to $\ep=1.5$.  Data sites for the RBF model are the ME points, while the data sites for the Fourier model are the MD points.} 
\label{fig:geom3d}
\end{figure}

As in the 2D results, we first examine the errors in reconstructing
the objects using the two parametric models.  Figure \ref{fig:geom3d}
displays the errors in reconstructing the objects as a function of the
square root of the number of data sites $N$. We use $\sqrt{N}$ since
these are 2D objects and thus the reciprocal of this value
gives a good measure of the spacing between data sites.  These errors
give a indication of the modeling capability of the RBF and Fourier
methods.  The results are similar to what we observed in 2D.  For the
smooth Object 1 (left plot in Figure \ref{fig:geom3d}), the RBF and
Fourier models are converging at a spectral rate, but at a much slower
rate for non-smooth Object 2.  The RBF and Fourier models are giving
similar errors for Object 1, with a few values of $N$ where the spherical
harmonic method is clearly better.  For Object 2, the RBF model
consistently gives better results than the spherical harmonic model as
$N$ increases. No direct comparison with the piecewise linear model is
given as the piecewise linear IB method always samples at the
interpolating points.

\subsubsection{Shape Parameter Study}
\label{sec:results3d_shape}

\begin{figure}[htbp]
\begin{center}
\includegraphics[width=3.0in]{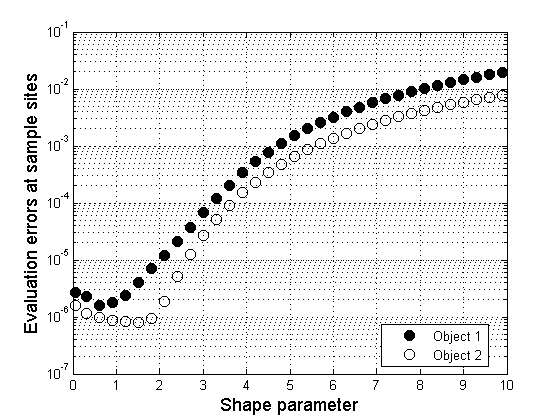}
\includegraphics[width=3.0in]{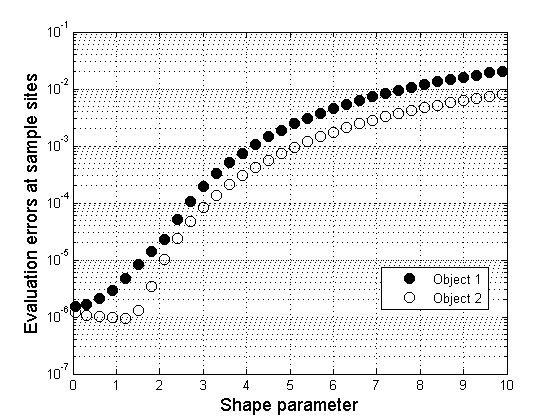}
\end{center}
\caption{Error in the shape at $1024$ sample sites
as a function of the shape parameter for $256$ data sites
on Object 1 (solid circles) and Object 2 (open circles) using
minimal energy points (left) and maximal determinant points (right)
for the data sites.} 
\label{fig:3dshape1}
\end{figure}

\begin{figure}[htbp]
\begin{center}
\includegraphics[width=3.0in]{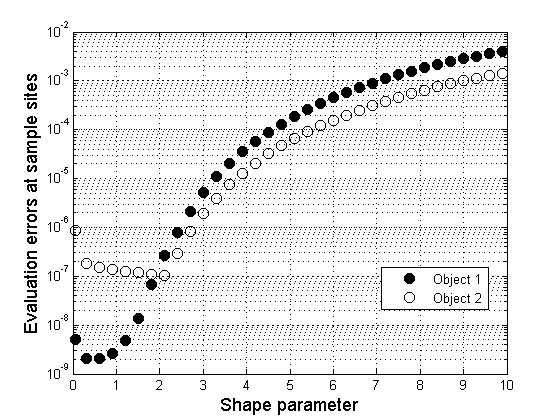}
\includegraphics[width=3.0in]{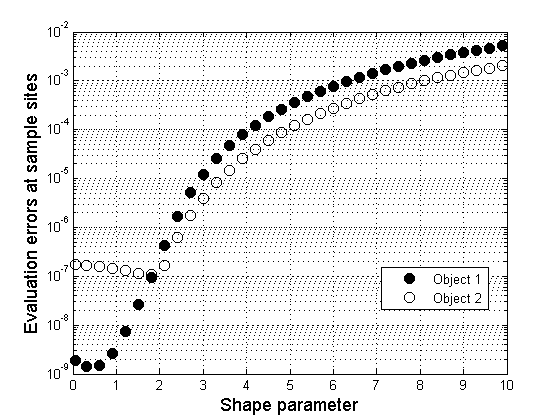}
\end{center}
\caption{Error in the shape at $1024$ sample sites
as a function of the shape parameter for $529$ data sites
on Object 1 (solid circles) and Object 2 (open circles) using
minimal energy points (left) and maximal determinant points (right)
for the data sites.} 
\label{fig:3dshape2}
\end{figure}

Figures \ref{fig:3dshape1} and \ref{fig:3dshape2} display the reconstruction errors of the RBF model for the two objects as a function of the shape parameter using $N=256$ data sites and $N=529$ data sites, respectively.  The left plot of each of these figures contains the results for the ME points, while the right plot contains the results for the MD points.  For $\ep \lesssim 0.85$, it was necessary to use the RBF-QR algorithm~\cite{FornbergPiret:2007} to compute the model in a numerically stable manner for the $N=529$ case.  We see similar results to the 2D shape parameter study from Section \ref{sec:results2d_shape}.  For the smooth Object 1 and the MD points the errors decrease rapidly as $\ep$ decreases and reach a minimum near $\ep=0$ (at $\ep=0$ in the $N=256$ case), which correspond to a spherical harmonic interpolant on these nodes.  For the ME points we see the error rise right as $\ep$ gets to zero.  For the non-smooth Object 2 and both types of nodes, we see that the error reaches a minimum at a larger value of $\ep$ that is well within the numerically safe range of RBF-Direct. The errors then increase slightly as $\ep$ decreases toward zero (with a jump up at $\ep=0$ in the case of the ME points).  From both Figures \ref{fig:3dshape1} and \ref{fig:3dshape2}, we see that the errors in the RBF model are much better for the MD points when $\ep$ is near zero, but as $\ep$ increases away from zero the errors are better for the ME points.

We make similar comments to those at the end of Section \ref{sec:results2d_shape} in regards to selection of the shape parameter for the 3D case, and thus refer the reader there.

\subsection{Comparison of Normal Vectors and Forces}
\label{sec:results3d_force}
We next focus on the errors in the parametric models in the approximation of the normal vectors to the objects and the forces.  In the case of computing the normal vectors, we compare the results to the traditional piecewise linear models based on triangulations of the surface.  As discussed in Section \ref{sec:implementation_pwl}, a comparison against the traditional piecewise linear 3D force model is not appropriate since this model is described purely algorithmically, and hence the underlying material constitutive model is not known and cannot be computed exactly even though the shape is known analytically.

\begin{figure}[htbp]
\begin{center}
\includegraphics[width=3.0in]{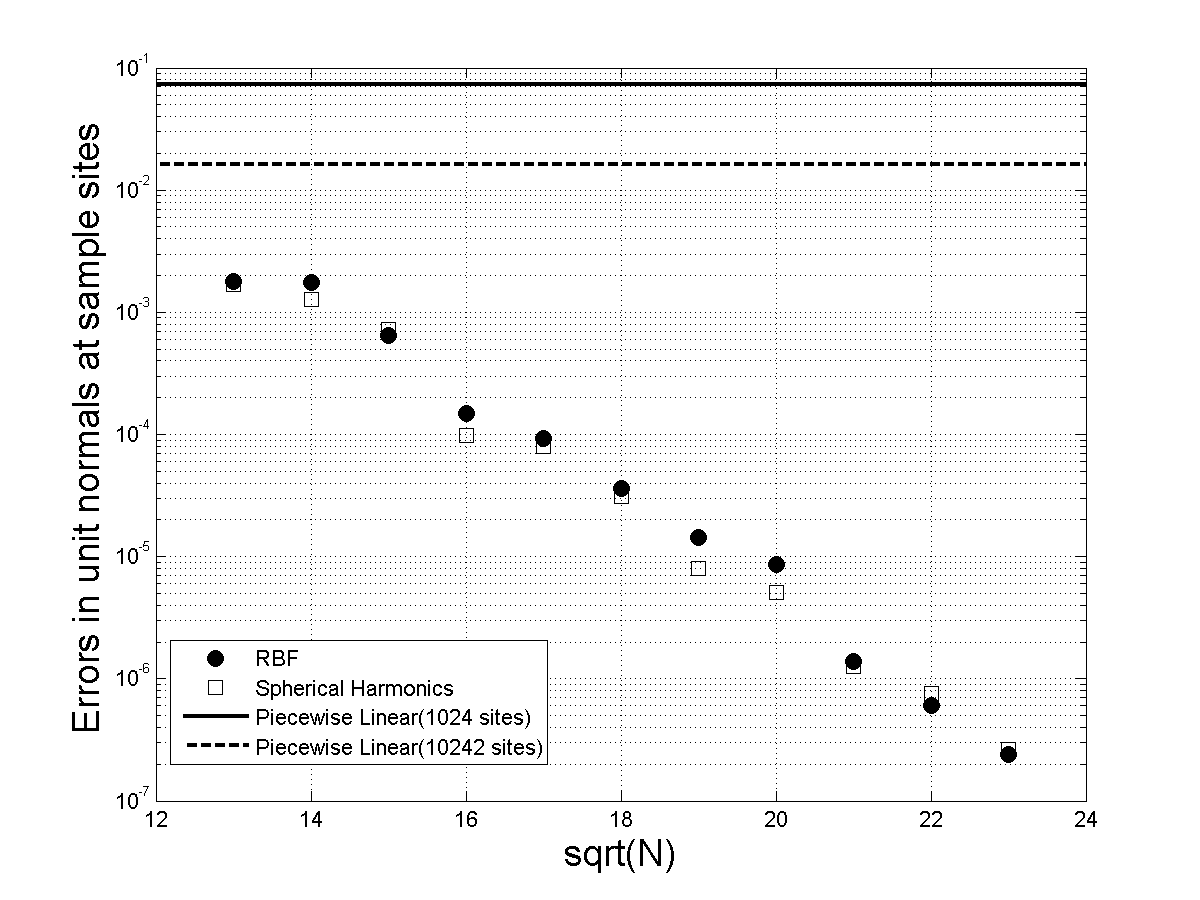}
\includegraphics[width=3.0in]{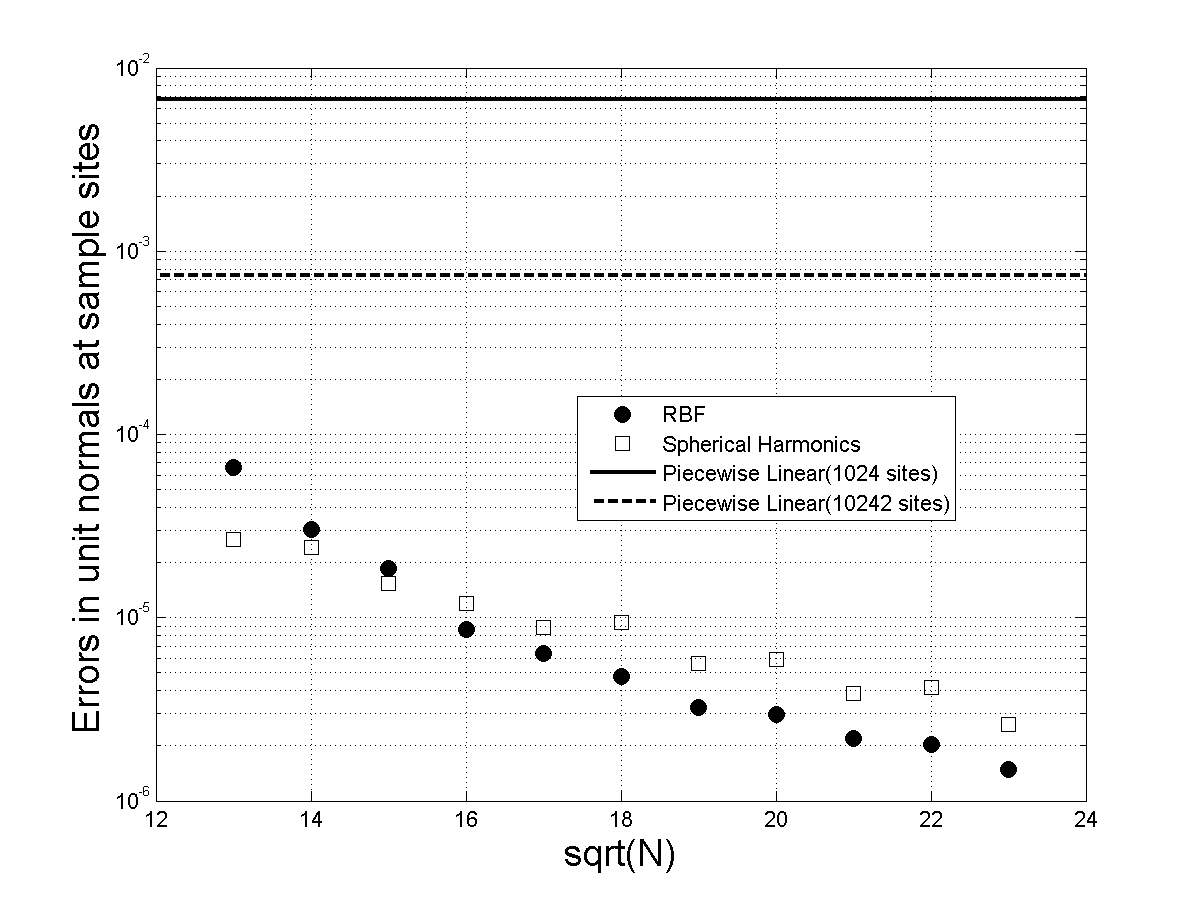} 
\end{center}
\caption{Errors in the approximations of the normal vectors to the 3D objects at $M=1024$ sample sites as a function of the square root of the number of data sites $N$.  The left plot is for Object 1, while the right one is for Object 2. The solid line denotes the error for the method used in the piecewise linear model with $1024$ IB points and the dashed line corresponds to the error with $10242$ IB points. Circles denote the errors for the RBF model and squares denote the Fourier model.  For the RBF model, $\ep=0.9$ for Object 1 and $\ep=1.5$. Data sites for the RBF model are the ME points, while the data sites for the Fourier model are the MD points.} 
\label{fig:normals3d}
\end{figure}

Figure \ref{fig:normals3d} displays the errors in the normal vectors at $M=1024$
sample sites as a function of the square root of the number of data sites $N$.  The solid and dashed lines in both plots from this figure denote the errors in the normal vectors at $1024$ and $10242$ IB points.  We see that increasing the number of IB points, decreases the errors in the normal vectors.  However, unlike the 2D case, both parametric models always give better results in the normal vector computations even for the high value of $10242$ IB points.  Additionally, the errors in these computations for Object 1 are similar for the RBF and Fourier models.  For Object 2, the RBF model gives consistently better results for increasing data sites $N$.

\begin{figure}[htbp]
\begin{center}
\includegraphics[width=3.0in]{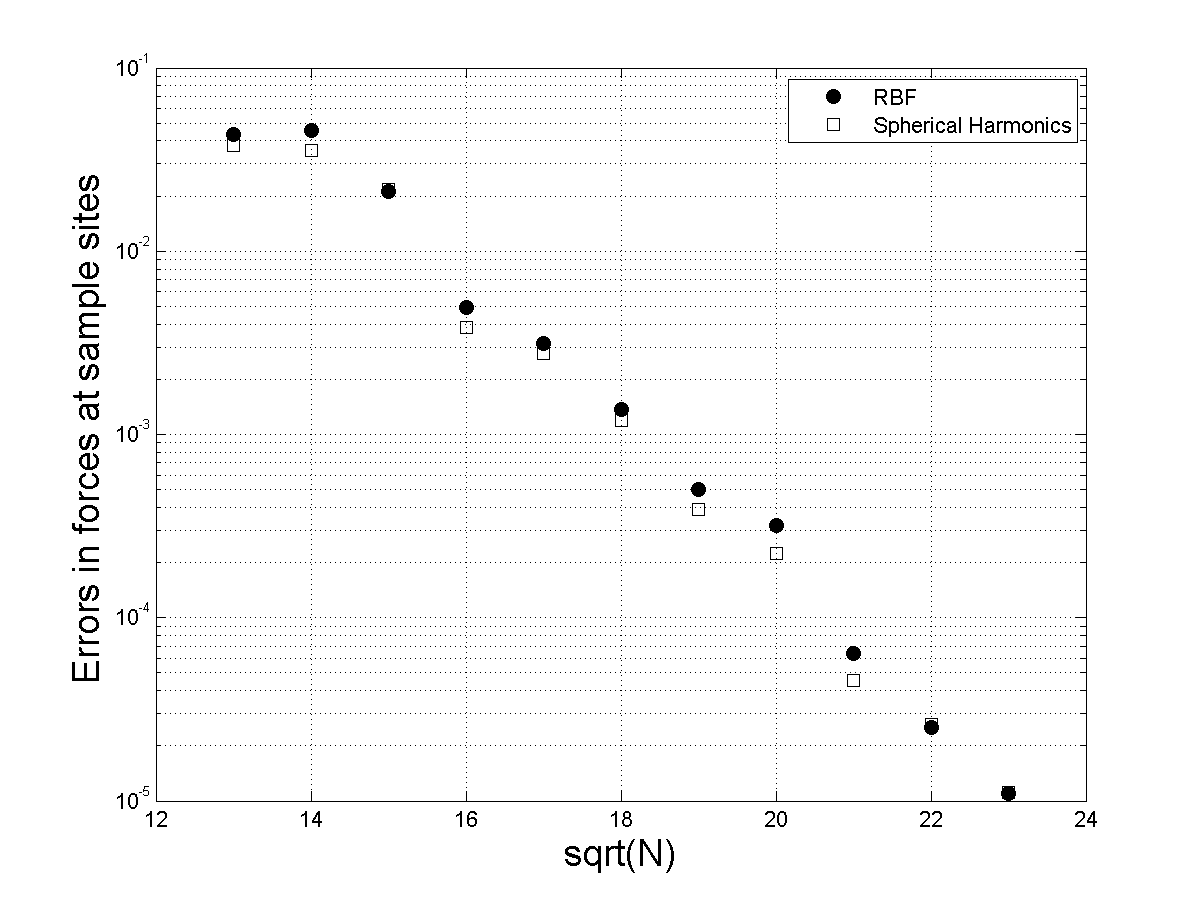}
\includegraphics[width=3.0in]{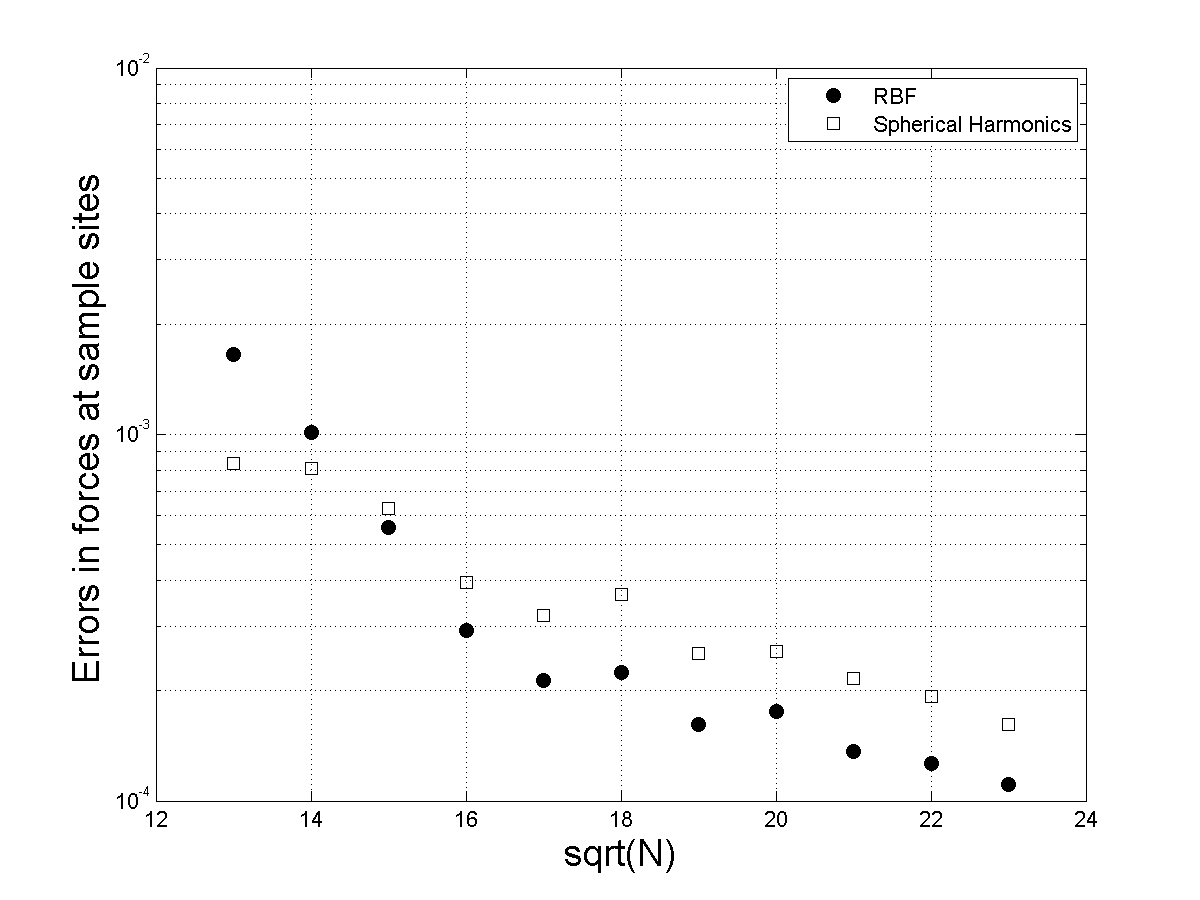} 
\end{center}
\caption{Errors in the approximations of the forces evaluated at $M=1024$ sample sites
as a function of the square root of the number of data sites $N$ 
for Object 1 (left) and Object 2 (right). Circles denote the errors for the RBF model and 
squares denote the Fourier model.  For the RBF model, $\ep=0.9$ for Object 1 and $\ep=1.5$. Data sites for the RBF model are the ME points, while the data sites for the Fourier model are the MD points.} 
\label{fig:force3d}
\end{figure}

We lastly focus on the errors in the computation of the forces that occur in both parametric models. Figure \ref{fig:force3d} displays the errors in forces evaluated at $1024$ sample sites as a function of the square root of the number of data sites $N$.  In all experiments, both the coefficient of surface tension $\gamma$ and the spring constant $K_0$ are set to 0.2.  We see that the results between smooth and non-smooth objects are consistent with those from the shape reconstruction and normal vector approximations.

\subsection{Comparison of the Computational Cost}
\label{sec:results3d_compcost}

We conclude the 3D results experiments by examining the computational cost associated with the three methods.  As in the 2D experiments, we measure the computational cost as the elapsed wallclock time required to compute the interpolation coefficients, evaluate the interpolants, compute the normal vectors and compute the forces.  We pre-compute and store the $LU$ decomposition of the spherical harmonic interpolation matrix \eqref{eq:sph_linsys} and the Cholesky decomposition $LL^T$ of the RBF interpolation matrix \eqref{eq:rbf_linsys}.  We also pre-compute matrices for evaluating the interpolants, the derivatives, and the force operator once the interpolation coefficients have been determined (see Section \ref{sec:implementation_param} for details). We do not account for these pre-computations in our timing results.  As in 2D, all computations were performed in \matlab using the machine described in Section \ref{sec:results2d_compcost}.

Since for the piecewise linear model the number of evaluation sites is the 
same as the number of data sites, the total computational cost includes only the time required to compute the normal vectors and forces (see Section \ref{sec:implementation_pwl} for details), we do not include the time to compute the triangulation of the surface. 

\begin{figure}[htbp]
\begin{center}
\includegraphics[width=3.0in]{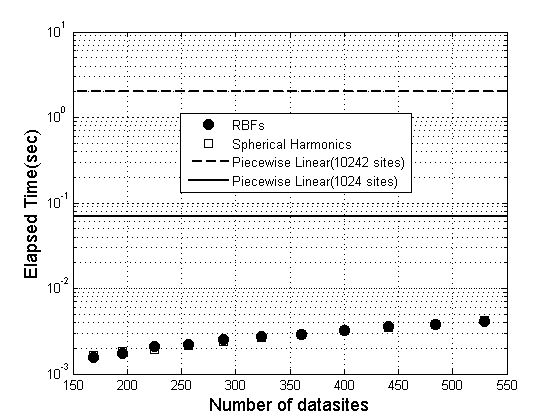}
\end{center}
\caption{Elapsed wallclock time (in seconds) for one object
to perform interpolation, evaluation, the computation of normal vectors, and 
the computation of forces at $M=1024$ sample sites as a function of the number 
of data sites $N$.  The piecewise linear computations were done with $1024$ IB points (solid line) and $10242$ IB points (dashed line) for comparison.} 
\label{fig:timings3d}
\end{figure}

Figure \ref{fig:timings3d} displays the elapsed time between the RBF, Fourier, and traditional piecewise linear models. The results for the RBF and Fourier models are displayed as a function of the number of data sites $N$ for a fixed number of $M=1024$ sample sites.  Two results are presented for the piecewise linear model: one with $1024$ IB points (solid) line and one with $10242$ IB points (dashed line).  We can see from the figure that the parametric models require significantly less time than the piecewise linear model, especially for the $10242$ case. For $N=529$ data sites, the parametric models are over one order of magnitude faster than the piecewise linear model with $1024$ IB points and nearly 3 orders of magnitude better with $10242$ IB points.

\section{Summary}\label{sec:summary}

The IB method is a common numerical methodology for applications involving fluid structure interactions.  Our particular interest in this method is in simulating platelet aggregation during blood clotting.  In this application, the platelets are modeled as immersed elastic structures whose shape changes dynamically in response to blood flow and chemistry. One of the fundamental ingredients of this application of the IB method (and many others involving immersed structures) is how to model the platelets geometrically, so that internal structural forces can be computed at specified locations on the platelet surface.  The current strategy is to use piecewise linear models for representing the platelets.  In this paper we have presented two alternative geometric models for platelets: RBFs and Fourier-based methods.  Both of these models are based on a parametric representation of the surface using polar coordinates in 2D and spherical coordinates in 3D.  This choice of parameterization is motivated by the observed shape of platelets both during their inactive and active states.  We have described how these new models can be used for constructing and maintaining the platelet's representation, computing the normal vectors to the platelet surface, and computing the internal structural forces.  We have presented numerical comparisons between the traditional piecewise linear models and the new RBF and Fourier-based models in both 2D and 3D.  Our findings indicate that both the RBF and Fourier methods provide viable alternatives to the traditional approach in terms of geometric modeling accuracy, force accuracy, and computational efficiency. 

Although both the RBF and Fourier-based methods provided comparable results in terms of error characteristics and computational efficiency, we would advocate the use of the RBF-based models for the following reasons:
\begin{itemize}
\item they are easier to implement; 
\item they have accuracy similar to that of Fourier methods for smoothly-perturbed objects with similar computational costs; 
\item they are more accurate than Fourier methods for roughly perturbed objects with similar computational costs;
\item they are more flexible than Fourier methods in terms of changing the underlying parameterizations of the objects (\emph{e.g.} changing to an elliptical parameterization rather than polar)~\cite{FuselierWright:2010}.
\end{itemize}
One issue with the RBF models is how to choose an appropriate shape parameter.  We will study this issue as part of our next step in applying the RBF-based models in an IB simulation.  This step will involve implementing the RBF-based models in a full IB simulation of platelet aggregation.  The simulation will require projection of the forces from the sample points to the Eulerian mesh, computation of the Navier-Stokes system with forcing based upon the platelets, and then movement of the platelets via updating of the RBF data points. We will study how the shape parameter affects the simulations and compare the results of these simulations to those based on the traditional piecewise linear models for platelets.

\vspace{0.2in} {\bf Acknowledgments:} We would like to acknowledge
useful discussions concerning this work within the CLOT group at the
University of Utah and with Prof. Robert Guy (UC-Davis).  
Support for this work came in part from NIGMS
grant R01-GM090203 (VS, ALF, RMK), from NSF grant DMS-0540779 (ALF,
GBW), and from NSF grant DMS-0934581 (GBW).
\bibliographystyle{elsart-num-sort}
\bibliography{Article}

\end{document}